\documentclass[12pt,leqno]{article}
\usepackage[latin1]{inputenc}
\usepackage{amsxtra}
\usepackage{amsmath}
\usepackage{amssymb}
\usepackage{amsfonts}
\usepackage[all]{xy}
\usepackage{mathrsfs}
\usepackage{amsthm}
\usepackage{enumerate}
\usepackage{mathabx}
\usepackage{graphicx}
\usepackage{mathtools}
 \usepackage[margin=1.5in]{geometry}
 \usepackage{tikz}
  \usepackage{epigraph}

\newtheorem{thm}[equation]{Theorem}
\newtheorem{cor}[equation]{Corollary}
   \newtheorem{Defi}[equation]{Definition}
\newtheorem{lem}[equation]{Lemma}
\newtheorem{prop}[equation]{Proposition}

\newtheoremstyle{example}{\topsep}{\topsep}%
     {}
     {}
     {\bfseries}
     {.}
     {2pt}
     {\thmname{#1}\thmnumber{ #2}\thmnote{ #3}}

   \theoremstyle{example}
   
   \newtheorem{rem}[equation]{Remark}

  \numberwithin{equation}{section}

\def\eps{{\varepsilon}}

\def\CC{\mathbb{C}}

\def\PP{\mathbb{P}}

\def\ZZ{\mathbb{Z}}

\def\gen{\mathfrak{g}}

\def\len{\mathfrak{l}}

\def\Cc{\mathcal{C}}

\def\Fc{\mathcal{F}}
\def\Gc{\mathcal{G}}

\def\Nc{\mathcal{N}}
\def\Hc{\mathcal{H}}
\def\Oc{\mathcal{O}}

\def\rb{\mathbf{r}}
\def\Db{\mathbf{D}}

\def\sb{\mathbf{s}}

 \def\Bl{{\on{Bl}}}
 
 \def\cHom{{\Hc om}}
 \def\Coker{\operatorname{Coker}\nolimits}
 \def\Com{\on{Com}}
 \def\cCom{\Cc om}

  \def\del{\partial}

 \def\End{\operatorname{End}\nolimits}

  \def\gCoh{{\on{gCoh}}}
 
 \def\gMod{{\on{gMod}}}
 \def\gVect{\on{gVect}}

  \def\Hom{\operatorname{Hom}\nolimits}
  \def\Hot{\Hc ot}

 \def\Id{\operatorname{Id}\nolimits}
 \def\Im{\operatorname{Im}\nolimits}

 \def\k{\mathbf k}
 \def\Ker{\operatorname{Ker}\nolimits}

 \def\Mat{\operatorname{Mat}\nolimits}
 \def\max{{\on{max}}}
 
\def\Mor{\operatorname{Mor}\nolimits}

\def\on{\operatorname}
\def\ol{\overline}

\def\NC{\on{NC}}

\def\Proj{{\on{Proj}}}
\def\PSS{{{\PP\SS}}}

\def\rk{\operatorname{rk}\nolimits}

 \def\Spec {\on{Spec}}
 \def\SS{{\on{SS}}}
 \def\St{{\on{St}}}

\def\ul{\underline}
\def\uHom{\ul{\on{Hom}}}

\def\wt{\widetilde}

\def\sb{\textbf{s}}
\def\tb{\textbf{t}}

\def\lra{\longrightarrow}

\def\(({(\hskip -1mm (}
\def\)){)\hskip -1mm )}
 \def\be{\begin{equation}}
\def\ee{\end{equation}}
\def\ed{\end{document}}

\def\1{{\mathbf {1}}}
\def\0{{\mathbf {0}}}

\title{ Complete complexes and spectral sequences}
\author{  Mikhail Kapranov, Evangelos Routis }
 
\begin{document}

\maketitle

\thanks{\em To Yuri Ivanovich Manin }

\begin{abstract}
By analogy with the classical (Chasles-Schubert-Semple-Tyrell)
 spaces of complete quadrics and complete collineations, we introduce  the variety
of complete complexes. Its points can be seen as 
equivalence classes of spectral sequences of a certain type. 
We prove that the set of such equivalence classes has a structure
of a smooth projective variety.  We show that it provides
a  desingularization, with normal crossings boundary, 
of  the Buchsbaum-Eisenbud variety of complexes, i.e., 
 a compactification of the union of its maximal strata. 

\end{abstract}

  \tableofcontents
  
  \vfill\eject

  \addtocounter{section}{-1}
 
 \section{Introduction}

 \noindent {\bf A. Background and motivation}.
 The spaces of complete collineations and complete quadrics form a beautiful and very important chapter
 of  algebraic geometry, going back to the classical works of Chasles and Schubert in the 19th century, see 
   \cite{semple}, \cite{tyrrell}, \cite{laksov:origin}, \cite{laksov:linear}, \cite{vainsecher}, \cite{thaddeus}, \cite{deconcini-macpherson},  \cite{kausz}
 and references therein. 
 They provide explicit examples of wonderful compactifications (i.e., of smooth compactifications
 with normal crossings boundary). 
 
 To recall the basic example, the group $\PP GL_n(\CC)= GL_n(\CC)/(\CC^*\cdot{\bf 1})$ has an obvious
  compactification by the projective space $\PP( \Mat_n(\CC))$ but it is not wonderful since the complement,
  the determinantal variety, is highly singular. Now let $V,W$ be two $\CC$-vector spaces of the same dimension
  $n$. A {\em complete collineation} from $V$ to $W$ is a sequence of the following data
  (assumed to be nonzero and considered each up to a non-zero scalar factor): 
  
  \begin{enumerate}
  \item[(0)] A linear operator $A=A_0: V\to W$, possibly degenerate (not an isomorphism). Note that
  $\Ker(A)$ and $\Coker(A)$ have the same dimension.
  
  \item[(1)] A linear operator $A_1: \Ker(A_0)\to\Coker(A_0)$, possibly degenerate.
  
  \item[(2)] A linear operator $A_2: \Ker(A_1)\to\Coker (A_1)$, possibly degenerate, and so on, until we obtain a non-degenerate linear operator. \end{enumerate}    
  One of the main results of the classical
   theory says that the set of  complete collineations has a natural structure of a smooth projective
   variety over $\CC$, containing $\PP GL_n(\CC)$ as an open subset ($A_0$ non-degenerate)
   so that the complement is a divisor with normal crossings. 
   
   \vskip .2cm
   
   We now want to look at this classical construction from a more modern perspective. We can view
   a linear operator $A: V\to W$ as a 2-term cochain complex. Then the sequence $(A_\nu)$ as above
   is nothing but a {\em spectral sequence}: a sequence of complexes $(E_\nu^\bullet, D^\nu)$
   such that each $E^\bullet_{\nu+1}$ is identified with the cohomology $H^\bullet_{D^\nu}(E_\nu^\bullet)$. 
   
   \vskip .2cm
   
  This suggests a generalization of the construction of complete collineations involving  more
  full-fledged (simply graded) spectral sequences. In this paper we develop such a generalization.
  The role of the group $GL_n(\CC)$ (or its projectivization $\PP GL_n(\CC)$) is played
  by appropriate strata in the 
    {\em Buchsbaum-Eisenbud variety of complexes}
  $C(V^\bullet)$ and its projectivization $\PP C(V^\bullet)$. Here $V^\bullet$ is a graded vector space and  $C(V^\bullet)$ consists
  of all ways of making $V^\bullet$ into a cochain complex, see \S \ref {subsec:affine} and 
 \cite{kempf}, \cite{deconcini-strickland} for more
  background. The varieties $C(V^\bullet)$ are known to share many important properties
  of determinantal varieties, in particular, they are {\em spherical varieties}: the action of the
  group $GL(V^\bullet) = \prod GL(V^i)$ on the coordinate ring has simple spectrum, i.e.,  each irreducible representation
  enters at most once.  
  
  \vskip .2cm
  
 \noindent{\bf B. Summary of results.}
 Our results can be summarized as follows. For simplicity, consider the projective variety of
 complexes $\PP C(V^\bullet)$. Let $\PP C^\circ(V^\bullet)$ be the union of its maximal 
 $GL(V^\bullet)$-orbits,
 a smooth open dense subvariety in $\PP C(V^\bullet)$, 
 see \eqref{eq:generic-str}.     
 
 At the same time let $\PP SS(V^\bullet)$ be the set of
 equivalence classes of  spectral sequences $(E_\nu^\bullet, D^\nu)$, $\nu=0,\cdots, N$, 
 of variable (finite)  length $N$, see \S  \ref{sec:CCSS}A, 
 in which:
 \begin{itemize}
  \item $E_0^\bullet=V^\bullet$.
  
  \item Each $D^\nu$, $\nu=0,\cdots, N-1$,  is not entirely zero and considered up to
 an overall scalar.
 
 \item The ``abutment" $E^\bullet_N$  
  does not admit any two consecutive nonzero spaces (so the spectral sequence must degenerate
 at $E_N^\bullet$). 
 
 \end{itemize}
 Then:
 \begin{enumerate}
 \item[(1)] The set $\PP SS(V^\bullet)$ admits the structure of 
 a smooth projective variety $\ol{\PP C}(V^\bullet)$ over $\CC$.
 
 \item[(2)]  $\ol{\PP C}(V^\bullet)$
  contains $\PP C^\circ(V)$ as an open dense part, and the complement
  $\PP SS(V^\bullet)-\PP C^\circ(V)$ is a divisor with normal crossings.
  
  \item[(3)] One can obtain $\ol{\PP C}(V^\bullet)$ as the successive blowup of 
  the closures of the natural strata in $\PP C(V^\bullet)$. 
 
 \end{enumerate}
  These results   are obtained by combining Theorems \ref{thm:wonder1} and \ref{thm:CCSS}.
  The realization of $\ol{\PP C}(V^\bullet)$ as an iterated blowup generalizes the approach of
  Vainsencher \cite{vainsecher} to complete collineations. In the main body of the paper we work over any algebraically closed field $\k$ of characteristic $0$ and consider the  varieties $C(V^\bullet)$ as well. Also, more generally, for any graded locally free sheaf of finite rank $V^\bullet$ over an arbitrary normal variety $X$ over $\k$, we introduce relative versions of varieties of complexes $ C_X(V^\bullet)$ and $\PP C_X(V^\bullet)$ (cf. Section \ref{subsec:relcom}). When $X$ is smooth we obtain the analogs of the results (1)-(3) above. 
  \vskip .2cm
  
  \paragraph{ C. Phenomena behind the results.} 
  The main phenomenon that makes our theory work, is the remarkable {\em self-similarity of the variety of
  complexes}. More precisely, $C(V^\bullet)$ is subdivided into strata (loci of complexes with prescribed ranks of the differentials).
  The transverse slice to a stratum passing through a point $D\in C(V^\bullet)$ (i.e., a differential in $V^\bullet$),
  is itself a variety of complexes but corresponding to the graded vector space $H^\bullet_D(V^\bullet)$
  of cohomology of $D$ (cf. Propositions \ref{prop:ncone} and \ref{prop:rel-ncone}). This generalizes the familiar self-similarity of the determinantal varieties: the transverse slice to the stratum formed by matrices of fixed rank inside a determinantal variety, is
  itself a determinantal variety of smaller size. In particular, our analysis implies 
    that our stratification is conical in the
  sense of MacPherson and Procesi \cite{macpherson-procesi}. 
  
  \vskip .2cm
  
  Further, the classical intuitive reason behind the appearance of  complete collineations has a transparent
   homological meaning.
  To recall this reason, consider a 1-parameter family $A(t)$ of linear operators $V\to W$ (depending,
  say, analytically on a complex number $t$ near $0$). If $A(t)$ is nondegenerate for $t\neq 0$
  but $A_0=A(0)$ is degenerate, then the ``next Taylor coefficient" of $A(t)$ gives $A_1: \Ker(A_0)\to\Coker
  (A_0)$,  the further Taylor coefficients give $A_2$ and so on. This gives the limit $\lim_{t\to 0} A(t)$
  in the space of complete collineations in the classical theory. 
  
  \vskip .2cm
  
  If  we now have an analytic 1-parameter family $D(t)$ of differentials in the same graded $\CC$-vector
  space $V^\bullet$, we can view the Taylor expansion of $D(t)$ as a  single differential  $\Db$ in the
  graded $\CC((t))$-vector space $V^\bullet\otimes_\CC \CC((t))$. The fact that $D(t)$ is analytic at $0$
  (so we are talking about Taylor, not Laurent expansions), means that $\Db$ preserves the
  $t$-adic filtration in  $V^\bullet\otimes_\CC \CC((t))$. The associated spectral sequence
  of the filtered complex
  is essentially simply graded, and it represents the limit of $(V^\bullet, D(t))$, as $t\to 0$, in
  our compactification. 
  
  \vskip .2cm

  \vskip .2cm
  
\noindent{\bf D. Future directions.} We expect our varieties of complete complexes to have interesting
enumerative invariants, generalizing the many remarkable properties of complete collineations.

Historically, the first example of a ``complete" variety
of geometric objects was the Chasles-Schubert space $\ol Q_n$ of complete quadrics, which gives a wonderful compactification of the variety $Q_n$
 of smooth quadric hypersurfaces in $\PP^n$, see \cite{deconcini-macpherson}.  
 From our point of view, $Q_n$ can be seen as a particular case of the {\em variety of self-dual complexes}.
 That is, we start with a
 graded (by $\ZZ$ or $\ZZ+{1\over 2}$) vector space $V^\bullet$  which is identified with its graded dual by 
 a graded symmetric bilinear form
 and consider  all  ways of making $V^\bullet$ into a self-dual complex. The corresponding analog
 of $\ol Q_n$ is then formed by the variety of {\em self-dual spectral sequences}. We leave its study to a future work.

    \vskip .2cm
    
  \noindent{\bf Acknowledgements.} This work was supported by the World Premier International Research Center Initiative (WPI Initiative), MEXT, Japan.

 
 \section{Categories of complexes}\label{sec:categories}

 Let $\k$ be an algebraically closed field and $\Lambda$ be a
 finitely generated  commutative $\k$-algebra. 
  We denote by $\gMod_\Lambda$ the category of finitely generated
  graded graded $\Lambda$-modules $V^\bullet$.
  That is,  $V^\bullet =\bigoplus_{i\in \ZZ}  V^i$, with all  $V^i$ finitely generated,
   and $V^i=0$ for $|i|\gg 0$.  For $j\in \ZZ$
 the shifted graded module $V^\bullet[j]$ is defined by 
  $(V^\bullet[j])^i = V^{j+i}$. 
 
 \vskip .2cm

 We denote by $\cCom_\Lambda$ the category of
 cochain complexes over $\Lambda$, i.e., of graded modules
  $V^\bullet=\bigoplus V^i\in\gMod_\Lambda$
 equipped with a differential $D$,  a collection of  $\Lambda$-linear maps $D_i: V^i\to V^{i+1}$
 satisfying $D_{i+1}\circ D_i=0$. We will consider $D=(D_i)$ as a morphism
 $V^\bullet \to V^\bullet[1]$ in $\gMod_\Lambda$. For a complex $(V^\bullet, D)$ we have the graded
 module $H^\bullet_D(V)$ of cohomology.  

 \vskip .2cm
 
 We define the shifted complex $(V^\bullet, D)[j]$ to have
 the underlying graded module $V^\bullet[j]$ as above and the differential 
 $D^{V^\bullet[j]}$ having components
 given by
 \be\label{eq:complex-shift}
 D^{V^\bullet[j]}_i \,\,=\,\, (-1)^{j} D_{j+i}^{V^\bullet}. 
 \ee
 One way to explain this formula is to represent $V^\bullet[j] = \Lambda[j]\otimes_\Lambda V^\bullet$
(here $\Lambda[j]$ is the ring $\Lambda$ put in degree $(-j)$). Then \eqref{eq:complex-shift}
corresponds to defining the differential via the graded Leibniz rule.

 \vskip .2cm
 
 We recall that a morphism  of complexes $f: (V^\bullet, D) \to (V'{}^\bullet,
 D')$ is called {\em null-homotopic}, if it is of the form $f= D's+sD$, where $s: V^\bullet\to V'{}^\bullet[-1]$
 is any morphism in $\gMod_\Lambda$. In this case we write $f\sim 0$. 
 For two morphisms of complexes
 $f,g: (V^\bullet, D) \to (V'{}^\bullet,
 D')$ we say that $f$ {\em is homotopic to}  $g$ and write $f\sim g$, if $f-g\sim 0$. 
 Null-homotopic morphisms form an ideal in $\Mor(\cCom_\Lambda)$,
 and the quotient category is called the {\em homotopy category of complexes} and will
 be denoted $\Hot_\Lambda$.
 
 \vskip .2cm
 
 \begin{Defi}\label{def:adm}
 (a) A complex $(V^\bullet, D)\in \cCom_\Lambda$ will be called {\em admissible}, if:
 \begin{enumerate}
 \item[(a1)]  Each $V^i$ is a projective $\Lambda$-module.
 
 \item[(a2)] Each image $\Im(D_i)\subset V^{i+1}$ is, locally, on the Zariski topology of $\Spec(\Lambda)$,
 a direct summand in $V^{i+1}$. 
   \end{enumerate}
   
   (b) Let $X$ be a $\k$-scheme of finite type. A complex $V^\bullet$ of coherent sheaves on $X$ will be
   called admissible, if, for any affine open $U\subset X$, the complex $\Gamma(U, V^\bullet)$ of
   modules over $\Lambda=\Oc(U)$, is admissible. 
 \end{Defi}

 \begin{prop}\label{prop:indec1}
 (a) For an admissible complex $(V^\bullet, D)$, each $H^i_D(V^\bullet)$ is a projective
 $\Lambda$-module. 
 
 \vskip .2cm
 
 (b) Assume that $\Lambda$ is a local ring. Then any admissible complex
  $V^\bullet\in\cCom_\Lambda$
 can be written as a direct sum $V^\bullet = A^\bullet
  \oplus H^\bullet$, where $A^\bullet$ is an admissible acyclic complex ($H^\bullet(A^\bullet)=0$),
  and $H^\bullet$ is an admissible  complex with zero differential
   (so $H^\bullet\simeq H^\bullet(V^\bullet)$).
   
   \vskip .2cm 
   
    (c) Further, for a local $\Lambda$,
     any acyclic admissible complex $A^\bullet$ is contractible (i.e., $\Id_{A^\bullet}\sim 0$). Therefore
     for any two admissible complexes $(V^\bullet, D)$ and $(V'{}^\bullet, D')$ over $\Lambda$ we have
  \[
 \Hom_{\Hot_\Lambda}\bigl( (V^\bullet,D),  (V'{}^\bullet, D')\bigr) \,\,\simeq 
 \,\, \Hom_{\gMod_\Lambda}(H^\bullet_D(V^\bullet), H^\bullet_{D'}(V'{}^\bullet)\bigr)
 \]
 \end{prop}
 
 \noindent {\sl Proof:} (a) Since $\Im(D_i)$ is locally a direct summand of a projective $\Lambda$-module
 $V^{i+1}$, it is projective. Therefore $\Ker(D_i)$ which is the kernel of the surjective
 morphism $D_i: V^i\to \Im(D_i)$ is itself projective and, moreover, locally a direct
 summand in $V^i$. So $\Im(D_{i-1})$ being, locally,  a direct summand in
 $V^i$, is in fact locally a direct summand in $\Ker(D_i)$, and so $H^i_D(V)=
 \Ker(D_i)/\Im(D_{i-1})$ is projective. 
 
 \vskip .2cm
 
 (b) Since $\Lambda$ is local, all locally direct summands discussed above are in fact
 direct summands of $\Lambda$-modules. Let $H^i\subset\Ker(D_i)$ be a direct
 complement to $\Im(D_{i-1})$, and let $W^i\subset V^i$ be a direct complement
 to $\Ker(D_i)$, so that $V^i=\Im(D_{i-1})\oplus H^i \oplus W^i$.  
  Then $H^\bullet\subset V^\bullet$ is a subcomplex
 with zero differential. Putting $A^i=\Im(D_{i-1})\oplus W^i$, we get a subcomplex
 $A^\bullet\subset V^\bullet$ which is acyclic, and $V^\bullet=A^\bullet\oplus H^\bullet$. 
 
 \vskip .2cm
 
 (c) If $A^\bullet$ is admissible and acyclic, we write, as before, $A^i=\Im(D_{i-1})\oplus W^i$.
 Since $\Im(D_{i-1})=\Ker(D_i)$, the restricted morphism ${D_i}_{|W^i}$ is an isomorphism $W^i\to\Im(D_i)$.
 Denote $s_{i+1}: A^{i+1}\to A^i$ to be the composite map
 \[
 A^{i+1}\buildrel \on{pr}\over\lra \Im(D_i) \buildrel {D_i}^{-1}_{|W^i}\over\lra W^i\hookrightarrow A^i.
 \]
 Then $s=(s_i: A^i\to A^{i-1})$ satisfies $Ds+sD=\Id_{A^\bullet}$. \qed
 
 \vskip .2cm
 
 When $\Lambda=\k$, all complexes are admissible, and we obtain:
 
 \begin{prop}\label{prop:indec2}
 Indecomposable objects in the abelian category $\cCom_\k$ are the following:
 \begin{enumerate}
 \item[(1)] $\k[j]$, $j\in\ZZ$;
 
 \item[(2)] $\{ \k\buildrel \Id\over\to\k\}[j]$, $j\in\ZZ$ (the 2-term complex with differential being the identity).
  \end{enumerate}
  \end{prop}
  
  \noindent {\sl Proof:}  A complex with trivial differential is a direct sum of summands of type (1).
  Further, an acyclic complex $A^\bullet$ is a direct sum of summands of type (2). Indeed, in the notation
  of the proof of Proposition \ref{prop:indec1}(c),     $A^\bullet$ splits into a direct
  sum of 2-term complexes $\{W^i\buildrel {D_i}_{|W^i}\over\longrightarrow\Im(D_i)\}$, which are thus direct sums of
  summands of type (2). So our statement follows from Proposition  \ref{prop:indec1}(b).   \qed


 \section{The affine variety of complexes, its strata and normal cones}   \label{subsec:affine}

 From now on we assume that the characteristic of $\k$ is $0$. 
 Let $V^\bullet=\bigoplus V^i$ be a finite-dimensional graded $\k$-vector
 space. 
 \begin{Defi}\label{defi:affine}The affine variety of complexes associated to $V^\bullet$  is the closed subscheme 
   \[
  C(V^\bullet) =\Big\{D=(D_i)\in\prod_i \Hom(V^i, V^{i+1})\Big{|} D_{i+1}\circ D_i=0 \,\,\text{for all}\,\,\, i\Big\}
  \]
 in the affine space $\prod_i \Hom(V^i, V^{i+1})$.
 \end{Defi}

 In other words $C(V^\bullet) $ is the subscheme of $\prod_i \Hom(V^i, V^{i+1})$ formed by all ways of making $V^\bullet$ into a complex. 
 See \cite {deconcini-strickland} for background. In particular, we note:
 
  \begin{prop} \label{normal}
  $C(V^\bullet)$ is a reduced scheme (affine algebraic variety). Further,
 each irreducible component of $C(V^\bullet)$ is a normal variety. \qed
 \end{prop}

 The group $GL(V^\bullet) = \prod_i GL(V^i)$ acts naturally on 
 $C(V^\bullet)$. An element $g=(g_i)$,
 $g_i\in GL(V^i)$, sends $D=(D_i)$ to
 \be\label{eq:gr-action}
 (gD)_i = g_{i+1}D_i g_i^{-1}.
 \ee
 Orbits of $GL(V^\bullet)$ on $C(V^\bullet)$ are nothing but isomorphism classes of complexes
  $(V^\bullet, D)$
 with all possible $D$. We will call these orbits the {\em strata} of $C(V^\bullet)$ and denote by $[D]$ the stratum
 passing through $D$. By the Krull-Schmidt theorem, strata (isomorphism classes)
  are labelled by the multiplicities of the
 indecomposable summands of $(V^\bullet, D)$ in the category $\cCom_\k$. Using the description
 of indecomposables given by Proposition \ref{prop:indec2}, one obtains an explicit
 combinatorial description of the strata. Let us recall this description, together with some further
 properties of strata and their closures that have been established in \cite {gonciulea}.
 
  \vskip .2cm
 
 Without changing the essense of the problem, we can (and will) assume that 
 $V^\bullet = \bigoplus_{i=0}^m V^i$ is concentrated in degrees $[0,m]$, and denote $n_i=\dim(V^i)$.
 Let $R$ be the set of sequences ${\rb} = (r_1, \cdots, r_m)$, $r_i\in\ZZ_{\geq 0}$, satisfying the
 conditions
 \[
 r_i + r_{i-1} \leq n_i, \quad i=0, \cdots, m+1.
 \]
 Here we put $r_0=r_{m+1} = n_{m+1}=0$. The set $R$ is partially ordered by
 \[
 {\rb} \leq {\rb}' \quad \Leftrightarrow \quad r_i\leq r'_i, \,\,\forall i. 
 \]
 For any $\rb\in R$ we denote
 \[
 \begin{gathered}
 C^\circ_\rb(V^\bullet) \,\,=\,\,\bigl\{ D\in C(V^\bullet)\bigr| \,
 \rk(D_i) = r_{i+1}, \,\, i=1, \cdots, m \bigr\},
 \\
  C_\rb(V^\bullet) \,\,=\,\,\bigl\{ D\in C(V^\bullet)\bigr| \,
 \rk(D_i) \leq r_{i+1}, \,\, i=1, \cdots, m \bigr\}. 
 \end{gathered}
 \]
 
  \begin{prop}[\cite{gonciulea}] \label{strata}
 \begin{itemize}
 \item[(a)]  The strata of $C(V^\bullet)$ are precisely the   $C^\circ_\rb (V^\bullet)$,
 $\rb\in R$. They are  non-empty, locally closed, smooth subvarieties. 
 
 \item[(b)] The closure of  $C^\circ_\rb (V^\bullet)$ is  $C_\rb (V^\bullet)$. In particular, each
  $C_\rb(V^\bullet)$ is irreducible.
  
  \item[(c)] We have  $C_\rb (V^\bullet)\subset  C_{\rb'} (V^\bullet)$,  if and only if
  $\rb\leq\rb'$.
  
  \item[(d)] The irreducible components of $C(V^\bullet)$ are precisely the $C_\rb (V^\bullet)$
  where $\rb$ runs over maximal elements of the poset $R$.   \qed
  
  \end{itemize}
  \end{prop}
  
  \begin{prop}\label{prop:sch-strata-abs}
  (a) The subvariety $C_\rb(V^\bullet)$ coincides with the subscheme in $C(V^\bullet)$
  given by the vanishing of the minors of size $r_{i+1}\times r_{i+1}$ of the differentials $D_i$
  for all $i$. In other words, the subscheme thus defined, is reduced. 
  
  \vskip .2cm
  
  (b) The scheme-theoretic intersection $C_\rb(V^\bullet) \cap C_\sb(V^\bullet)$
  coincides with the variety $C_{\min(\rb, \sb)}(V^\bullet)$, where
  \[
  \min(\rb, \sb) \,\,=\,\,(\min(r_1, s_1), \cdots, \min(r_m, s_m)). 
  \]
   \end{prop}
  
  \noindent{\sl Proof:} Part (a) is one of the main results of De Concini-Strickland
   \cite {deconcini-strickland}. Part (b) follows from the following well known  property of the determinantal
   ideals in the ring $\k[a_{ij}]$ of polynomials in the entries of an indeterminate $p\times q$ matrix
   $\|a_{ij}\|$. The ideal generated by all $r\times r$ minors contains the ideal generated by
   all $(r+1)\times (r+1)$ minors. \qed

  \vskip .2cm
  
    Let $Y$ be a closed subscheme of a $\k$-scheme $Z$ of finite type, and $I_Y\subset\Oc_Z$
  the sheaf of ideals of $Y$. We denote by
  \[
  \Nc^*_{Y/Z} = I_Y/I_Y^2, \quad \Nc_{Y/X} = \uHom_{\Oc_Y}(\Nc_{Y/Z}, \Oc_Y)
  \]
  the conormal and normal sheaves to $Y$ in $Z$. We will be particularly interested in the
  case when $\Nc^*_{Y/Z}$ (and therefore $\Nc_{Y/Z}$) is locally free, i.e., represents
  a vector bundle on $Y$. The total space of this vector bundle is then a scheme which we
  call the {\em normal bundle} to $Y$ in $Z$ and
  denote
  \[
  N_{Y/Z}\,\,=\,\,\on{Spec}_{\Oc_Y} S^\bullet_{\Oc_Y}  (I_Y/I_Y^2). 
  \]
We further denote by
\[
\NC_{Y/Z} \,\,=\,\, \on{Spec}_{\Oc_Y} \biggl( \bigoplus_{n=0}^\infty I_Y^n/I_Y^{n+1}\biggr)
\]  
 the {\em normal cone} to $Y$ in $Z$. Because of the surjection of sheaves of $\Oc_Y$-algebras
 \[
  S^\bullet_{\Oc_Y}  (I_Y/I_Y^2) \lra \bigoplus_{n=0}^\infty I_Y^n/I_Y^{n+1},
 \] 
 $\NC_{Y/Z}$ is a closed subscheme in $N_{Y/Z}$. In particular, $\NC_{Y/Z}$ is a
 ``cone bundle" over $Y$: it is equipped with an affine morprhism $\NC_{Y/Z}\to Y$
 whose fiber over a $\k$-point $y\in Y$ is a cone $(\NC_{Y/Z})_y$ in the linear
 space $(N_{Y/Z})_y$. 
  
 The above constructions extend easily to the case when $Y$ is locally closed (instead
 of closed) subscheme in $Z$. In this case, we define
 \[
\Nc^*_{Y/Z}=\Nc^*_{Y/Z^\circ},\quad N_{Y/Z}= N_{Y/Z^\circ}\quad \text{and} \quad \NC_{Y/Z}=\NC_{Y/Z^\circ},
 \]
 where $Z^\circ\subset Z$ is any open subset containing $Y$ and such that $Y$ is closed in $Z^\circ$.
 See \cite{fulton}  for more background.
 
  \vskip .2cm
 
 We now specialize to the case when $Z=C(V^\bullet)$ and $Y=[D]$ is the stratum
 through a $\k$-point $D$. Let ${\rb} = (r_1, \cdots, r_m)\in R$ and let
  $*$ stand for any of the categories $\cCom_\k, \Hot_\k$ or $\gMod_\k$. 
We  denote by 
 \[
 \Hom_{*}^{\leq\rb}\bigl( (V^\bullet, D), (V^\bullet, D)[1]\bigr)\,\,\subset \,\,
 \Hom_{*}\bigl( (V^\bullet, D), (V^\bullet, D)[1]\bigr)
 \]
  the closed subvariety    formed by    morphisms $f=(f_i:V^i\rightarrow V^{i+1})$ 
  such that  $\rk(f_i) \leq r_i$ for all $i$. 
  
   \begin{prop}\label{prop:TDCom}
   
   \begin{enumerate}
\item[(a)]  The Zariski tangent space to $C(V^\bullet)$ at $D$ is found by:
 \[
 T_D C(V^\bullet) \,\,=\,\, \Hom_{\cCom_\k}\bigl( (V^\bullet, D), (V^\bullet, D)[1]\bigr). 
 \]
 \item[(b)]  Suppose $D$ is contained in some $C_\rb (V^\bullet)$. Then, the Zariski tangent 
 space to $C_\rb(V^\bullet)$ at $D$ is found by:
  \[
 T_DC_\rb (V^\bullet) \,\,=\,\, \Hom_{\cCom_\k}^{\leq\rb}\bigl( (V^\bullet, D), (V^\bullet, D)[1]\bigr). 
 \]
 \end{enumerate}
 \end{prop}
 
 \noindent{\sl Proof:} (a)  By definition, $T_D C(V^\bullet)$ is the set of points $D_\eps$ of
 $C(V^\bullet)$ with values in $\k[\eps]/\eps^2$ which extend the $\k$-point $D$. 
 Since $C(V^\bullet)$ is embedded into the affine space $\Hom_{\gVect_\k}(V^\bullet, V^\bullet[1])$, 
 we can write $D_\eps = D+\eps f$ where $f\in \Hom_{\gVect_\k}(V^\bullet, V^\bullet[1])$.
 The condition for $D_\eps$ to be a point of $C(V^\bullet)$ is the vanishing of
 $D_\eps^2 = D^2 + \eps(Df+fD)$. Since $D^2=0$ by assumption, we are left with
 $Df+fD=0$ which, in virtue of the convention \eqref{eq:complex-shift}, means that 
 $f: (V^\bullet, D)\to (V^\bullet, D)[1]$ is a morphism of complexes. Part (b) is similar. 
  \qed
  
   \begin{prop}\label{prop:TDD}
   The Zariski tangent space $T_D[D] \subset T_D C(V^\bullet)$ consists of those
 morphisms of complexes $(V^\bullet, D)\to (V^\bullet, D)[1]$, which are homotopic to $0$. 
   \end{prop}
  
  \noindent {\sl Proof:}  By definition, $[D]=GL(V^\bullet)\cdot D$ is the orbit of $D$ under the
  action \eqref{eq:gr-action}. Therefore
  \[
  T_D[D] \,\,=\,\,\Im\bigl\{{\mathfrak {gl}}(V^\bullet) \buildrel h_D\over\lra T_D C(V^\bullet)
  \bigr\}
  \]
 is the image of the Lie algebra ${\mathfrak {gl}}(V^\bullet)$ under the infinitesimal action $h_D$
 induced by  \eqref{eq:gr-action}. To differentiate  \eqref{eq:gr-action}, we take
 \[
 g=(g_i), \,\,\, g_i = 1+\eps s_i, \,\,\, s_i \in {\mathfrak gl}(V^i), \,\,\, \eps^2=0. 
 \]
Then
\[
(gD)_i \,\,=\,\, (1+\eps s_{i+1}) D_i (1-\eps s_i) \,\,=\,\,
D_i + \eps (s_{i+1} D_i - D_i s_i)
\] 
 which is precisely a perturbation of $D$ by a morphism homotopic to $0$. Since $s_i$ can
 be arbitrary, the statement follows.
   \qed

 \vskip .2cm
 
 Since $[D]$ is an orbit of $GL(V^\bullet)$, the conormal sheaf $\Nc^*_{[D]/ C(V^\bullet)}$
 is locally free and so we can speak about the normal bundle $N_{[D]/C(V^\bullet)}$.
 Proposition \ref{prop:indec1}(c)  together with the above implies:
 
 \begin{cor}
 The fiber at $D$ of the normal bundle $N_{[D]/C(V^\bullet)}$ is found by:
 \[
 (N_{[D]/C(V^\bullet)})_D \,\,\simeq \,\, \Hom_{\gMod_\k}\bigl( H^\bullet_D(V^\bullet),
 H^\bullet_D(V^\bullet)[1]\bigr). \qed
 \]
  \end{cor}
  
  \begin{prop} \label{prop:ncone}
  \begin{enumerate}
  \item[(a)] The fiber of $\NC_{[D]/C(V^\bullet)}$ at $D$ is found by:
  \[
 (\NC_{[D]/C(V^\bullet)})_D \,\,\simeq \,\,
  C(H^\bullet_D(V^\bullet)) \,\,\subset \,\,\Hom_{\gMod_\k} \bigl( H^\bullet_D(V^\bullet), 
  H^\bullet_D(V^\bullet)
  [1]\bigr).
  \]
  
  \item[(b)]  Suppose $D$ is contained in some $C_\rb (V^\bullet)$. Then, the fiber of
   $\NC_{[D]/C_\rb (V^\bullet)}$ at $D$ is found by:
   \[
 (\NC_{[D]/ C_\rb (V^\bullet)})_D \,\,\simeq \,\,
  C_{\rb-\rb_D} (H^\bullet_D(V^\bullet)) \,\,\subset \,\,\Hom_{\gMod_\k}^{\leq\rb-\rb_D} 
  \bigl( H^\bullet_D(V^\bullet), H^\bullet_D(V^\bullet)
  [1]\bigr),
  \]
  where $\rb_D$ is the sequence of ranks of the differential $D$.
  \end{enumerate}
  \end{prop}

  \noindent {\sl Proof:} (a) We note first that
  \be\label{eq:cone-factor}
   (\NC_{[D]/ C(V^\bullet)})_D \,\,\simeq \,\, \NC_{D/C(V^\bullet)}\bigl/ T_D [D],
  \ee
  the quotient of the tangent cone at the point $D$ by the action of the vector space $T_D[D]$
  (we consider this space as an algebraic group acting on the tangent cone).
  Now, the normal cone to any subvariety $Y$ in any $X$ at a $\k$-point $y$ is found, inside
  $T_yY$, by:
  
  \begin{itemize}
  \item [(1)]  Considering all $\k[[\eps]]$-points $x_\eps$ of $X$ which are (1st order) tangent to $Y$
  
  \item[(2)] Restricting the equations of $Y$ in $X$ to such $x_\eps$.
  
  \item[(3)] Equating to $0$ the next (after the linear) lowest nonvanishing terms in $\eps$ in these equations.
  \end{itemize}
  
  In our case $X=\Hom_{\gMod_\k}(V^\bullet, V^\bullet[1])$ is an affine space, so  in Step (1) above 
  it is enough to take the $\k[[\eps]]$-points of the form $D_\eps = D+\eps f$ with 
  $f\in \Hom_{\gVect_\k}(V^\bullet, V^\bullet[1])$. By Proposition \ref{prop:TDCom},
  for $D_\eps$ to be 1st order tangent to $C(V^\bullet)$, it is necessary and sufficient
  that $f$ be a morphism in $\cCom_\k$, not just in $\gMod_\k$.

  Further, in Step(2), the equations of $C(V^\bullet)$ after restricting
  to $D_\eps$ are the matrix elements of
  \[
  D_\eps^2 \,\,=\,\, D^2 + \eps (Df+fD) + \eps^2 f^2,
  \]
  so the next nonvanishing coefficient in Step (3) is $f^2$. This means that
  \[
  \NC_{D/C(V^\bullet)} \,\,=\,\,\bigl\{ f\in \Hom_{\cCom_\k} \bigl( (V^\bullet, D), (V^\bullet, D)[1]\bigr)
  \,\bigr| \, f^2=0\bigr\}.
  \]
  Our statement now follows from \eqref{eq:cone-factor}, the identification of $T_D[D]$
  in Proposition \ref{prop:TDD} and Proposition \ref{prop:indec1}(c).
  
  \vskip .2cm

(b) The isomorphism of Proposition \ref{prop:indec1}(c) restricts to an isomorphism of the subspace $
 \Hom_{\Hot_\k}^{\leq \rb }\bigl( (V^\bullet,D),  (V^\bullet, D)[1]\bigr) $ with
 $ \Hom_{\gMod_\k}^{\leq\rb-\rb_D}(H^\bullet_D(V^\bullet), H^\bullet_{D}(V^\bullet)[1]\bigr)$.
 So it suffices to repeat the argument of (a).\qed


  \section{The projective variety of complexes}\label{subsec:projective}
  
  It follows from Definition \ref{defi:affine}  that $C(V^\bullet)$ is given by homogeneous (quadratic) equations in the
  linear space $\Hom_{\gMod_\k}(V^\bullet, V^\bullet[1])$; those are the matrix elements of all the maps $D_{i+1}\circ D_i: V^i\to V^{i+2}$, where $D_i \in \Hom(V^i, V^{i+1})$ for all $i$. We can therefore give the following definition:
  
  \begin{Defi}  The  projective variety of complexes associated to $V^\bullet$ is the projectivization
  \[
  \PP C(V^\bullet) \,\,\subset \,\, \PP \bigl(\Hom_{\gMod_\k}(V^\bullet, V^\bullet[1])
  \bigr)
  \]
  of $C(V^\bullet)$.
\end{Defi}
  
  It follows that $\PP C(V^\bullet)$ is a reduced scheme (projective algebraic
  variety), and each of its irreducible components is normal.

   For any nonzero differential $D=(D_i: V^i\to V^{i+1})$ in 
    $C(V^\bullet)$, we  denote its image in $\PP C(V^\bullet)$ by $\PP D$. 
    Further, for any  ${\bf 0}\neq \rb\in R$, we denote by $\PP C^\circ_\rb(V^\bullet)$ and
    $\PP C_\rb(V^\bullet)$ the images of the stratum $C^\circ_\rb(V^\bullet)$
    and of its closure in $\PP C(V^\bullet)$ respectively. We call the 
     $\PP C^\circ_\rb(V^\bullet)$
    for $\rb\neq {\bf 0}$,   the {\em strata} of $\PP\Com(V^\bullet)$ 
    and denote by $[\PP D]$ the stratum
 passing through $\PP D$.
 
 The properties of the affine varieties of complexes and their strata imply at once:
 
 \begin{prop}\label{prop:proj-strata}
  \begin{itemize}
 \item[(a)]  The strata of  $\PP C(V^\bullet)$ are precisely the  $GL(V^\bullet)$-orbits.
  They are  non-empty, locally closed, smooth subvarieties. 
 
 \item[(b)] The closure of  $\PP C^\circ_\rb (V^\bullet)$ is  $\PP C_\rb(V^\bullet)$.
  In particular, each
  $\PP C_\rb (V^\bullet)$ is irreducible.
  
  \item[(c)] We have  $\PP C_\rb(V^\bullet)\subset  \PP C_{\rb'} ( V^\bullet)$ if and only if
  $\rb\leq\rb'$.
  
  \item[(d)] The irreducible components of $\PP C(V^\bullet)$ are precisely the 
  $\PP C_\rb (V^\bullet)$
  where $\rb$ runs over maximal elements of the poset $R$.   \qed
  
  \end{itemize}
 \end{prop}

    \begin{prop} \label{prop:pncone}
  \begin{enumerate}
  
  \item[(a)] The fiber of $\NC_{[\PP D]/\PP C(V^\bullet)}$ at $\PP D$ is found by:
  \[
 (\NC_{[\PP D]/\PP C(V^\bullet)})_{\PP D} \,\,\simeq \,\,
  C(H^\bullet_{ D}(V^\bullet)) .
  \]
  
  \item[(b)]  Suppose $\PP D$ is contained in some $\PP C_\rb (V^\bullet)$.
   The fiber of\ $\NC_{[\PP D]/\PP C_\rb (V^\bullet)}$ at $\PP D$ is found by:
   \[
 (\NC_{[\PP D]/ \PP C_\rb (V^\bullet)})_{\PP D} \,\,\simeq \,\,
  C_{\rb-\rb_D}(H^\bullet_{ D}(V^\bullet)). \qed
  \]
  \end{enumerate}
  \end{prop}


 \section{The relative  varieties of complexes and their normal cones}\label{subsec:relcom} 
 
 Let $X$ be a normal algebraic variety over $\k$. We denote by $\gCoh_X$ the category of
 graded coherent sheaves $\Fc^\bullet = \bigoplus_{i\in\ZZ} \Fc^i$ and morphisms preserving
 the grading. For any two $\Fc^\bullet, \Gc^\bullet\in\gCoh_X$ we have a coherent sheaf
 $\uHom_{\gCoh_X}(\Fc^\bullet,\Gc^\bullet)$ on $X$ (local homomorphisms). 
 
 \vskip .2cm

 Let $V^\bullet = \bigoplus_{i=0}^m V^i$ be a graded vector bundle
 (locally free sheaf of finite rank) on $X$ concentrated in degrees $[0,m]$. For any $\k$-point $x\in X$ we denote by
 $V^\bullet_x$ the graded vector space obtained as the fiber of $V^\bullet$ at $x$. 
 The constructions of 
 \S \ref{subsec:affine} admit obvious relative versions. 
 \begin{Defi} \label{defi:rel-var-com}
 The relative affine variety of complexes associated to $V^\bullet$ is the $X$-scheme $C_X(V^\bullet)\buildrel \pi\over\to X$ that
  represents the functor $(Sch/X)\to(Set)$ given by
 $$T \,\,\mapsto \,\,\bigl\{ (\phi, D)\, | \phi: T\to X \text{ a morphism of schemes, }
  D \text{ a differential in } \phi^*V^\bullet\bigr\}. $$
 \end{Defi}
By definition, $C_X(V^\bullet)$ carries the {\em universal complex} of vector bundles
 \[
 (\ul V^\bullet, \ul D), \quad \ul V^\bullet = \pi^* V^\bullet, \quad \ul D =
  \text{the universal differential coming from Definition \ref{defi:rel-var-com}}.
 \]
  In particular, the fiber of $\pi$ at a $\k$-point $x\in X$ is the variety of complexes
 $C(V^\bullet_x)$. 
 
 \vskip .2cm
 
 We denote by 
 \[
M  \,\,= \,\,  \cHom_{\gCoh_X} (V^\bullet, V^\bullet[1])\,\,=\,\, \Spec_{\Oc_X} S^\bullet\biggl( \bigoplus_{i=0}^{m-1} V^i\otimes (V^{i+1})^*\biggr)
 \]
  the total space of the vector bundle
 $\uHom_{\gCoh_X}(V^\bullet, V^\bullet[1])$  (``relative space of matrices"). Then $C_X(V^\bullet)$ is a closed  conic subvariety in $M$. Let 
 $\PP M\to X$ be the projectivization of $M$ over $X$. 
 \begin{Defi} The relative projective variety of complexes associated to $V^\bullet$ is the projectivization $\PP C_X(V^\bullet)\subset\PP M$ of $C_X(V^\bullet)$. 
 \end{Defi}
 
 \vskip .2cm
 
 As before, for any $\rb\in R$ we denote by $C^\circ_{X,\rb}(V^\bullet )$, (resp. by $C_{X, \rb}(V^\bullet)$) the     locally closed,
 (resp. closed)
 subvariety in $C_X(V^\bullet)$ formed by differentials  $D=(D_i)$  with  the rank of $D_i$ equal to $r_i$ everywhere
 (resp. $\leq r_i$ everywhere). We refer to the $C^\circ_{X,\rb}(V^\bullet )$ as the {\em strata} of $C_X(V^\bullet)$.
 
 If $\rb\neq {\bf 0}$, we denote by $\PP C_{X,\rb}^\circ(V^\bullet)$, resp. $\PP C_{X,\rb}(V^\bullet)$,  the image of
 $C^\circ_{X,\rb}(V^\bullet)$, resp. $C_{X, \rb}(V^\bullet)$ in $\PP C_X(V^\bullet)$. We call the $\PP C_{X,\rb}^\circ(V^\bullet)$
 the {\em strata} of $\PP C_X(V^\bullet)$. 
 
 \begin{prop}\label{prop:C_X-properties}
 (a) $C_X(V^\bullet)$ and $\PP C_X(V^\bullet)$ are reduced schemes; each irreducible component of these schemes
 is a normal variety. 
 
   \vskip .2cm
 
 (b) If $X$ is smooth, then each stratum $C_{X,\rb}^\circ(V^\bullet)$, $\PP C_{X,\rb}^\circ(V^\bullet)$ is smooth.  

  \vskip .2cm
 
  (c) The subvariety $C_{X,\rb}(V^\bullet)$ coincides with the subscheme in $C_X(V^\bullet)$
  given by the vanishing of the minors of size $r_{i+1}\times r_{i+1}$ of the differentials $D_i$
  for all $i$.  
  
  \vskip .2cm
  
  (d) The scheme-theoretic intersection $C_{X,\rb}(V^\bullet) \cap C_{X,\sb}(V^\bullet)$
  coincides with the variety $C_{X,\min(\rb, \sb)}(V^\bullet)$, where
  \[
  \min(\rb, \sb) \,\,=\,\,(\min(r_1, s_1), \cdots, \min(r_m, s_m)). 
  \]

    \end{prop}
    
  \noindent{\sl Proof:} Parts (a)-(c) follow from the absolute case,  cf. Proposition  \ref {prop:sch-strata-abs}.
  The proof of (d) is also similar to the proof of Proposition  \ref {prop:sch-strata-abs}(b).
  \qed
  
  \vskip .2cm
 
\noindent  Let $S$ be a stratum of $C_X(V^\bullet)$. We denote by ${\ul V^\bullet}_S = {\ul V^\bullet}|_S$ the restriction 
to  $S$
 of the universal complex $\ul V^\bullet$. By the definition of the strata,  ${\ul V^\bullet}_S$ is an admissible complex
 (Def. \ref{def:adm}), and therefore its graded sheaf of cohomology 
 with respect to
 the restriction of the differential $\ul D$,  is locally free (a graded vector bundle).
  We denote this graded vector bundle by
 $H^\bullet_S:= H^\bullet_{\ul D}({\ul V^\bullet}_S)$.
  Note that $H^\bullet_S$ descends canonically
 to a graded vector bundle on the projectivization $\PP S$, which we will  denote by $H^\bullet_{\PP S}$. 
 
 Propositions \ref {prop:ncone}
 and \ref{prop:pncone}, describing the  normal cones to the strata fiberwise, can be formulated in a neater, global way,
 using relative varieties of complexes. The proofs are identical and we omit them.
 
 \begin{prop}\label{prop:rel-ncone}
 (a) Let $S = C^\circ_{X, \rb'}(X)$ be a stratum in $C_X(V)$ and $\rb'\leq \rb$.  Then
 \[
 \NC_S(C_X(V^\bullet)) \,=\, C_S(H^\bullet_S), \quad \NC_S(C_{X, \rb}(V^\bullet)) \,=\, C_{S, \rb-\rb'}(H^\bullet_S).
 \]
 
\noindent  (b)  Let $\PP S = \PP C^\circ_{X, \rb'}(X)$ be a stratum in $\PP C_X(V)$ and $\rb'\leq \rb$.  Then
 \[
 \NC_{\PP S}(\PP C_X(V^\bullet)) \,=\, C_{\PP S}(H^\bullet_{\PP S}), \quad \NC_{\PP S}(\PP C_{X, \rb}(V^\bullet)) \,=\, 
 C_{{\PP S}, \rb-\rb'}(H^\bullet_{\PP S}). \qed
 \]
 \end{prop}

 
 \section{Charts in varieties of complexes}
 We keep the notation of \S \ref{subsec:relcom}. 
 The goal of this section is to prove:
 
  \begin{prop}\label{prop:chart} 
  \begin{enumerate}
 \item[(a)] Let $c$ be a $\k$-point of  $C_X(V^\bullet)$ (resp. $\PP C_X(V^\bullet)$)
  belonging to a stratum $S$.
  There exists an isomorphism of an \'etale neighborhood of $c$   in $C_X(V^\bullet)$ 
  (resp.  in $\PP C_X(V^\bullet$) with an \'etale neighborhood of $c$   in $\NC_{S/C_X(V^\bullet)}$
  (resp. in  $\NC_{S/\PP C_X(V^\bullet)}$).
  
  \item[(b)]  Suppose that $c$ (and therefore $S$) is contained in some $C_{X,\rb}(V^\bullet)$ 
  (resp. $\PP C_{X, \rb}(V^\bullet)$). Then the isomorphism of part (a) restricts to an isomorphism
   of an \'etale neighborhood of $c$ in $C_{X, \rb}(V^\bullet)$ (resp.  in $\PP C_{X, \rb}(V^\bullet)$) with an
    \'etale  neighborhood of $c$  in $\NC_{S/C_{X, \rb}(V^\bullet)}$(resp. $\NC_{S/\PP C_{X, \rb}(V^\bullet)}$).
  \end{enumerate}
 \end{prop}
 
 \begin{rem}
(a)    Combining Propositions \ref{prop:chart}  and  \ref{prop:rel-ncone}, we obtain the following conclusion
  (``self-similarity of varieties of complexes''). 
  Any 
 relative affine or projective variety of complexes is modeled, near any point $c$, by another affine variety
 of complexes, which is typically simpler than the original one (depending on the singular nature of $c$).

 \vskip .2cm
 
 (b) The proposition implies, in particular,
  that the natural stratifications of $C_X(V^\bullet)$ and $\PP C_X(V^\bullet)$
 are conical in the sense of MacPherson and Procesi \cite{macpherson-procesi}.

 \end{rem}

 \vskip .2cm
 
 \noindent {\sl Proof of Proposition  \ref{prop:chart}:} We begin by a series of reductions.
  First, we only need to prove the statement for the affine variety of complexes:
   the projective case follows immediately from that by descent. 
  
  Second, we  need only to give the proof for part (a):   part (b) will  follow by  identical arguments.

  Third, it is enough to consider only the absolute case when $X$ is a point. Indeed, 
  by restricting, if necessary,  to an open subset of $X$, we can  assume that
   the graded vector bundle $V^\bullet$ is trivial, identified with $X\times V^\bullet_x$,
   where $V^\bullet_x$ is the fiber at some $x\in X$. In this case $C_X(V^\bullet) = X\times C(V^\bullet_x)$,
   and the stratum $S$ has the form $X\times S_x$, where $S_x$ is a stratum in $C(V^\bullet_x)$.
   A point $c$ has then the form $c=(x, D)$ where $x\in X$ and $D$ lies in $S_x$. 
   A chart for $C_X(V^\bullet)$ near $c$ will follow from a chart for $C(V^\bullet_x)$ near $D$. 
   
   \vskip .2cm

   So we assume that $X=\Spec(\k)$ and $V^\bullet$ is a graded vector space.
   Our point $c$ is therefore just a differential $D$ in $V^\bullet$ and $S=[D]$.  The remainder of the
   proof is subdivided into three steps. 
   
      \vskip .2cm

 \textbf{Step 1:} We write $V^{\bullet}=  H^{\bullet} \oplus A^{\bullet} $, where
  $H^{\bullet}$ is a complex with zero differential and $A^{\bullet}$ is an acyclic complex 
  (see Proposition \ref{prop:indec1}). Let $D_A$ be the differential of $A^{\bullet}$. Then we can write $D$ in matrix form as
 \be\label{eq:D-A}
 D\,\,=\,\, 
\left(
\begin{array}{c|c}
0 & 0  \\
\hline
0 & D_A\\
\end{array}
\right)
\ee
where the zero upper left part corresponds to $H^\bullet$. 
Now, we have an embedding of $C(H^\bullet)$ into $C(V^\bullet)$ defined by
\be\label{eq:D-delta}
\delta
\,\,
\mapsto  \,\, D_\delta :=
\left(
\begin{array}{c|c}
\delta & 0  \\
\hline
0 & D_A\\
\end{array}
\right), \quad \delta\in  C(H^\bullet). 
\ee
 Further,  $H^\bullet$ is isomorphic to $H_D^\bullet(V^\bullet)$, the cohomology of $V^\bullet$.
 Therefore, by Proposition \ref{prop:ncone}, we deduce that $C(H^\bullet)$ is isomorphic to the fiber 
 of the normal cone $\NC_{[D]/C(V^\bullet)}$ over $D$.  
 In other words, we have embedded the fiber of the normal cone over $D$ back into $C(V^\bullet)$.
  Our goal in the steps to follow is to extend this embedding to an \'etale map from  an open
   neighborhood of the fiber to an open neighborhood of $D$ in $C(V^\bullet)$.
 
 \vskip .2cm
 
 \noindent \textbf{Step 2:} 
  Let $\St_D$ be the stabilizer of $D$ for the action of $GL(V^\bullet)$ on $C(V^\bullet)$. 
  We write $g\in GL(V^\bullet)$ as  $g=(g_i: V^i\to V^i)_{i=0}^m$, and then write
each $g_i$ in the matrix form with respect to $V^i=H^i\oplus A^i$:

\[
g_i=\left(
\begin{array}{c|c}
P_i & Q_i  \\
\hline
R_i & S_i\\
\end{array}
\right), \quad P_i: H^i\to H^i, \,\, R_i: H^i\to A^i, \text{ etc.}
\]
The condition that $g\in\St_D$ means, in virtue of \eqref{eq:D-A} and  the action law \eqref{eq:gr-action}:
\be\label{eq:stab}
Q_{i+1} D_i=D_i R_{i}=0,  \quad S_{i+1}D_i=D_iS_{i},
\ee
which is a set of linear homogeneous equations on the matrix elements of the $g_i$.
In other words, \eqref{eq:stab} defines a linear subspace
$W\subset \End(V^\bullet) = \bigoplus_i\End_\k(V^i)$, 
and $\St_D=W\cap GL(V^\bullet)$. Choose a complementary affine subspace $L$ to $W$
passing through ${\bf 1}\in\End(V^\bullet)$ (which is the unit element of $GL(V^\bullet)$).
Then set $U=L\cap GL(V^\bullet)$.

\begin{lem}\label{lem:I}
The action map $\phi: U\to [D]$, $g\mapsto g\cdot D$, is birational, and its differential at ${\bf 1}\in U$
is an isomorphism. 
\end{lem}
  
  \noindent {\sl Proof:} It is clear from the construction (the tangent space  $T_\1 U$ is a complement 
  to $T_\1 \St_D$),
   that $\dim(U)=\dim[D]$, and $d_\1\phi$ is an isomorphism.  To see that $\phi$ is birational, 
   fix a generic $g_0\in U$
   and see how many $g\in U$ are there such that $\phi(g)=\phi(g_0)$. The latter condition means
    $g\cdot D=g_0\cdot D$,
   i.e., $g_0^{-1} g\in\St_D$, or, in other words, $g\in g_0\cdot \St_D$. Since $\St_D$
    is the intersection of $GL(V^\bullet)$
   with a linear subspace in $\End(V^\bullet)$, the coset $g_0\cdot\St_D$ also has this property,
    so $L\cap (g_0\cdot \St_D)$
   typically consists of one point, i.e. $\phi$ is generically bijective onto its image. Since $[D]$ is normal and the ground field $\bf k$ has characteristic 0, Zariski's Main Theorem for quasifinite morphisms implies that $\phi$ is generically an open immersion, i.e. birational. \qed

  \vskip .2cm
  
  Let us now extend \eqref{eq:D-delta} to an embedding 
  \be\label{eq:eta}
  \eta: \Hom(H^\bullet, H^\bullet[1]) \hookrightarrow 
  \Hom(V^\bullet, V^\bullet[1]), \quad \delta\mapsto D_\delta
  \ee
  by allowing $\delta$ to be an arbitrary morphism of graded vector spaces
  $H^\bullet\to H^\bullet[1]$ (not necessarily
  satisfying $\delta^2=0$).
  
  \begin{lem}\label{lem:II}
  Inside the tangent space $T_D \Hom(V^\bullet, V^\bullet[1])$ we have $\Im(d_0\eta)\cap T_D[D]=0$. 
  \end{lem}
  
  \noindent {\sl Proof:} Identifying the tangent space in question with the
   vector space $\Hom(V^\bullet, V^\bullet[1])$, we have
$$
  \Im(d_0\eta) =\Hom(H^\bullet, H^\bullet[1]) \subset \Hom (H^\bullet\oplus A^\bullet,
  H^\bullet[1]\oplus A^\bullet[1]) =\Hom(V^\bullet, V^\bullet[1]). 
  $$
  Note that $\Im(d_0\eta)$ consists of morphisms of complexes, not just of graded vector spaces, 
  since $H^\bullet$
  has zero differential and is a direct summand in $V^\bullet$ as a complex. 
  On the other hand, $T_D[D]$
  consists, by Proposition \ref{prop:TDD}, of morphisms of complexes $V^\bullet\to V^\bullet[1]$ 
  which are
  homotopic to $0$. Such morphisms induce the zero map on the cohomology. On the other hand,
  morphisms from $\Im(d_0\eta)$ are faithfully represented by their action on the cohomology, since
  $H^\bullet \simeq H^\bullet_D(V^\bullet)$. Therefore the intersection of the two subspaces is $0$.
   \qed 
  
  \begin{lem}\label{lem:III}
  The action map
  \[
  \Psi: U\times\Hom(H^\bullet, H^\bullet[1]) \lra 
  \Hom(V^\bullet, V^\bullet[1]), \quad (g,\delta) \mapsto g\cdot D_\delta,
  \]
  is, at  the point $(\1, 0)$, \'etale onto its image.
  \end{lem} 
  
  \noindent {\sl Proof:} As both the source and target of $\Psi$ are smooth,
   it is enough to show that the differential
  $d_{(\1,0)}\Psi$ is an injective linear map.  Now, 
  \[
  T_{(\1,0)}(U\times\Hom(H^\bullet, H^\bullet[1])) \,\,=\,\,T_\1 U\oplus \Hom(H^\bullet, H^\bullet[1]).
  \]
   The restriction of $d_{(\1,0)}\Psi$ to  the   summand $T_\1U$  is the map  $d_\1\phi$ which, by 
    Lemma \ref{lem:I}, maps it isomorphically
    to $T_D[D]$. The restriction of $d_{(\1,0)}\Psi$ to the summand
     $\Hom(H^\bullet, H^\bullet[1])$ is the embedding
    $d_0\eta$, see \eqref {eq:eta}. So, by Lemma \ref{lem:II}, its image does not intersect 
    the image of the
    summand $T_\1 U$, which is $T_D [D]$. This means that the map from the direct sum
     of the two
    summands is injective. \qed
    
    \begin{cor}
    The action map
    \[
    \Phi: U\times C(H^\bullet) \lra C(V^\bullet), \quad (g,\delta) \mapsto g\cdot D_\delta,
    \]
    (the restriction of $\Psi$) is  \'etale at the point $(\1,0)$. 
    \end{cor}
    
  \noindent {\sl Proof:} By Lemma \ref{lem:III}, $\Phi$ is, at $(\1,0)$, \'etale onto its image. 
  This image is 
  contained in $C(V^\bullet)$. Now, we look at 
 the irreducible components  $K$ of $C(V^\bullet)$ through $D$.  Applying Proposition \ref{strata}(d),
 we see that they are in bijection with irreducible components of
 $C(H^\bullet)$ through $0$ and therefore with irreducible components $K'$ of
  $U\times C(H^\bullet)$ through $(\1,0)$.
 The dimension of each component $K'$ is equal to the dimension of the corresponding $K$, 
 and we see that
  $K$ is covered by $K'$ near $D$ (and therefore $\Phi(K')=K$).  So $\Im(\Phi)=C(V^\bullet)$.   \qed

  \vskip .2cm

\noindent\textbf{Step 3:} The isomorphism $H^\bullet\simeq H^\bullet_D(V^\bullet)$ induces,
by Proposition \ref{prop:ncone}, an identification
\[
\xi: C(H^\bullet) \to C(H^\bullet_D(V^\bullet)) \,\,=\,\, \bigl(\NC_{[D]/C(V^\bullet)}\bigr)_D. 
\]
We extend it to a map
\[
\Xi: U\times C(H^\bullet) \lra NC_{[D]/C(V^\bullet)}, \quad (g,\delta)\mapsto g\cdot \xi(D_\delta),
\]
using the action of $GL(V^\bullet)$ on $NC_{[D]/C(V^\bullet)}$. This morphism is birational and, moreover,
biregular near $\{\1\}\times C(H^\bullet)$ by Lemma \ref{lem:I}, since $g\in U$ takes the fiber of $NC_{[D]/C(V^\bullet)}$
over $D$ to the fiber over $g\cdot D$. 

Now, the composition $\Phi\circ\Xi^{-1}: NC_{[D]/C(V^\bullet)}\to C(V^\bullet)$ is a rational map,
 regular and \'etale
at  the point $D\in NC_{[D]/C(V^\bullet)}$. Proposition \ref{prop:chart} is proved.


   \section{Complete complexes via blowups}\label{sec:blowups}

Let $X$ be a smooth irreducible variety over $\k$, and $V^\bullet$ be a graded vector bundle
over $X$ with grading situated in degrees $[0,m]$, and $\on{rk}(V^i)=n_i$. 
 We keep all the other notations of  \S \ref{subsec:relcom}. 
 In this section we construct
 the {\em variety of
complete complexes}, resp. {\em projective variety of
complete complexes} 
   by successively blowing up closures of strata in $C_X(V^\bullet)$ and
$\PP C_X(V^\bullet)$ respectively. 

\vskip .2cm

\noindent {\bf A. Reminder on blowups.} Let $Z$ be a scheme of finite type over $\k$
and $Y\subset Z$ a closed subscheme with sheaf of ideal $I_Y$. The {\em blowup of $Y$
in $Z$} is the scheme  
\[
\Bl_Y(Z) \,\,=\,\, \Proj \biggl( \bigoplus_{n=0}^\infty I_Y^n\biggr). 
\]
See \cite{fulton} and \cite{hartshorne} for general background. In particular, we have a natural projection
\[
p: \Bl_Y(Z)\to Z, \quad p^{-1}(Y) \,\,=\,\,\PP \NC_{Y/Z} \,\,:=\,\, \Proj 
\biggl( \bigoplus_{n=0}^\infty I_Y^n/I_Y^{n+1}\biggr), 
\]
which restricts to an isomorphism $\Bl_Y(Z) - p^{-1}(Y)\to Z-Y$. 
 We will be especially interested in the case when $Y$ is a smooth algebraic variety and the conormal sheaf
 $\Nc^*_{Y/Z}$ is locally free, in which case $\PP\NC_{Y/Z} \subset \PP N_{Y/Z}$
 is a closed subscheme in the projectivization of the normal bundle $N_{Y/Z}$.   Compare with
 \S \ref{subsec:affine}. 
 
 \vskip .2cm
 
 If $W\subset Z$ is a closed  subscheme, 
 the {\em strict transform} of $W$ is defined as 
  $W^{\on{st}} = \Bl_{W\cap Y}(W)$ where $W\cap Y$ is the scheme-theoretic intersection.
   It is a closed subscheme in $\Bl_Y(Z)$. 
   
   If $Z$ is an algebraic variety and  $W$ is an irreducible subvariety, then
   $W^{\on{st}}$ is  equal to the closure in $\Bl_Y(Z)$
 of $p^{-1}(W-Y)$ (in particular, it is empty, if $W\subset Y$). 
 More generally, if $W\subset Z$ is any subvariety, its
 {\em dominant transform} $\wt W$ is defined as \cite{li-li}:
 \[
  \wt W \,\,=
 \begin{cases}
 W^{\on{st}}, & \text{if } W\not\subset Y;
\\
\text{total inverse image } p^{-1}(W), & \text{if }  W\subset Y. 
 \end{cases}
 \]
 
 \begin{prop}\label{prop:blow-disjoint}
 Let $Z$ be a scheme of finite type over $\k$ and  $W_1, W_2\subset Z$ be closed subschemes.
 Let $Y=W_1\cap W_2$ (scheme-theoretic intersection). Then the
 strict transforms of $W_1$ and $W_2$ in $\Bl_Y(Z)$ are disjoint. 
   \end{prop}
   
   \noindent{\sl Proof:} This is Exercise 7.12 in \cite{hartshorne}. \qed
 

 \paragraph  { B. Details on the poset of strata.}   Recall the poset $R$ of integer vectors
 $\rb=(r_1, \cdots, r_m)$ labelling the strata (as well as the closures of the strata)
  of $C_X(V^\bullet)$, see \S \ref{subsec:affine}
 for the absolute case,  extended in
 \S \ref{subsec:relcom} to the relative case. Thus the zero vector $\0$ is the minimal element 
 of $R$. 
 For $\rb\in R$ we denote $|\rb|=\sum r_i$, and
 call this number the {\em length} of $\rb$. We denote 
 $R_l = \{\rb\in R: \, |\rb|=l\}$, and similarly for $R_{\geq l}$, $R_{<l}$ etc. 
 
 \begin{prop}
 (a) The poset $R$ is ranked with rank function $|\rb|$. That is, for any $\rb\leq\rb'$
 all maximal chains of strict inequalities
 $\rb=\rb^{(0)} < \cdots < \rb^{(l)}=\rb'$ have the same cardinality $l+1$, where
 $l=|\rb'|-|\rb|$. 
 
 \vskip .2cm
 
 \noindent (b) The set of minimal elements of $R_{\geq l}$ coincides with $R_l$. 
 \end{prop}
 
 \noindent {\sl Proof:} Both statements follow from the next property which is obvious
 from the definition of $R$ by inequalities.
 
 \begin{lem}
 If $\rb\in R$ and $r_i\neq 0$, then $\rb-e_i\in R$,
 where $e_i$ is the $i$th basis vector ($1$ at the position $i$, zeroes everywhere else). \qed
 \end{lem}
 
 \vskip .2cm
 
 We now introduce the notation for some subsets of $R$:
 
 $R^\max$ denotes the set of maximal elements  of $R$
 (which label irreducible components of $C_X(V^\bullet)$ 
 as well as of  $\PP C_X(V^\bullet)$). 
 
 $R^\flat=R-R^\max$ (labels closures of non-maximal  strata in $C_X(V^\bullet)$,  to be blown up).
 
 $R^\flat_{ >l}= R^\flat \cap R_{>l}$ etc. In particular, $R^\flat_{>0}$ labels
closures of non-maximal strata in $\PP C_X(V^\bullet)$.

\vskip .2cm

\begin{Defi}A graded vector space $V^\bullet$ is {\em sparse}, if the numbers $n_i=\dim(V^i)$ 
are such that
$n_in_{i+1}=0$ for each $i$.
 \end{Defi}
 
\begin{rem} $V^\bullet$ is sparse, iff $C(V^\bullet)=\{0\}$. 
\end{rem}

\begin{prop}\label{prop:max-sparse}
A vector $\rb\in R$ lies in $R^\max$ if and only if for each $D\in C^\circ_\rb(V^\bullet)$
the graded vector space $H^\bullet_D(V^\bullet)$ is sparse. 
\end{prop}

\noindent{\sl Proof:} Indeed, the fiber of the normal cone to the stratum $[D] =  C^\circ_ {\rb}(V^\bullet)$
 over $D$ is, by Proposition \ref{prop:ncone}, identified with $C(H^\bullet_D(V^\bullet))$. Saying that
 $[D]$ is a maximal stratum is equivalent to saying that its normal cone consists of just the
 zero section. \qed
 
 \vskip .2cm
 
 We denote by
 \be\label{eq:generic-str} 
 C^\circ(V^\bullet) \,\,=\,\, \bigcup_{\rb\in R^\max} C^\circ_{\rb}(V), \quad 
 \PP C^\circ(V^\bullet) \,\,=\,\, \bigcup_{\rb\in R^\max} \PP C^\circ_{\rb}(V)
 \ee
 the union of the maximal strata.  It is a smooth open dense subvariety in $C(V^\bullet)$, resp.
 $\PP C(V^\bullet)$. We will refer to it as 
  the   {\em generic part } of $C(V^\bullet)$, resp.
 $\PP C(V^\bullet)$.

 \paragraph { C. Main constructions and results.} Let $d$ be the maximal value of $|\rb|$
 for $\rb\in R$. Our first result gives a construction of a series of blowups of the varieties
 of complexes with good properties (at each step we perform a blowup with a smooth center). 
 For convenience of inductive arguments, we formulate the results  for both affine and
 projective varieties of complexes. 
 
 \begin{thm}\label{thm:blow}
 There exist towers of  blowups
 \[
 \begin{gathered}
 \ol{C}_X(V^\bullet) = C_X^{(d)}(V^\bullet)\to \cdots\to C_X^{(1)}(V^\bullet) \to 
 C_X^{(0)} (V^\bullet) = C_X(V^\bullet),
 \\
 \ol{\PP C}_X(V^\bullet ) = {\PP C}^{(d)}_X(V^\bullet) \to\cdots\to
 \PP C_X^{(2)}(V^\bullet) \to  \PP C_X^{(1)}(V^\bullet) =  \PP C_X(V^\bullet)
 \end{gathered}
 \]
 with the following properties (which define them uniquely). 
 For $\rb\in R$   let 
 $C_{X,\rb}^{(l)}(V^\bullet)\subset C_X^{(l)}(V^\bullet)$ be the iterated dominant transform of
 $C_{X,\rb}(V^\bullet)\subset C_X(V^\bullet)$, and for $\rb\in R_{> 0}$ let
 $\PP C_{X,\rb}^{(l)}(V^\bullet)\subset \PP C_X^{(l)}(V^\bullet)$ be the iterated dominant 
 transform of
 $\PP C_{X,\rb}(V^\bullet)\subset \PP C_X(V^\bullet)$. 
 \begin{enumerate}
 
 \item[(a)]   For any given $l$ and  for $|\rb|=l$, 
 the subvarieties $C_{X,\rb}^{(l)}$ (resp. $\PP C_{X,\rb}^{(l)}$)
 are smooth and disjoint.
 
 \item[(b)] We have
 \[
  C_X^{(l+1)}(V^\bullet) \,=\,\Bl_{\coprod\limits_{|\rb|=l} C_{X,r}^{(l)}(V^\bullet)} C_X^{(l)}(V^\bullet), \quad 
\PP  C_X^{(l+1)}(V^\bullet) \,=\,\Bl_{\coprod\limits_{|\rb|=l} \PP C_{X,r}^{(l)}(V^\bullet) }
\PP C_X^{(l)}(V^\bullet). 
 \]
  \end{enumerate}
 \end{thm}
 
 The theorem will allow us to make the following definition.
  
 \begin{Defi}\label{defi:complete} The relative 
 variety of complete complexes associated to $V^\bullet$,  is the variety $\ol{C}_X(V^\bullet)$.
 The relative projective variety of complete complexes associated to $V^\bullet$, is the variety
 $\ol{\PP C}_X(V^\bullet)$.

 \end{Defi}

 Our second result says that the iterated blowup we construct, provides  wonderful compactifications
 of the open strata in the irreducible components of  the    varieties of complexes. 
 
 \begin{thm}\label{thm:wonder1}
 \begin{enumerate}
 \item[(a)]  The variety $\ol{C}_X(V^\bullet)$ (resp. $\ol{\PP C}_X(V^\bullet)$) is smooth and equal to the disjoint
 union of the $C_{X,\rb}^{(d)}(V^\bullet)$ (resp, of the $\PP C_{X, \rb}(V^\bullet)$) for
 $\rb\in R^\max$. 
 
 \item[(b)] The varieties $\Delta_\rb = C^{(d)}_{X,\rb}(V^\bullet), \rb\in R^\flat$, are smooth and
  form a divisor with normal
 crossings in $\ol{C}_X(V^\bullet)$ which we denote $\del\ol{C}_X(V^\bullet)$. 
 Similarly, the  varieties 
 $\PP\Delta_\rb=\PP C^{(d)}_{X,\rb}(V^\bullet), \rb\in R^\flat_{>0}$,  are smooth and form a divisor with normal
 crossings in $\ol{\PP C}_X(V^\bullet)$, which we denote $\del\ol{\PP C}_X(V^\bullet)$,
 
 \item[(c)] The complement $\ol C_X(V^\bullet) - \del\ol C_X(V^\bullet)$
 is identified with the  disjoint union of the strata $C^\circ_{X,\rb}(V^\bullet)$ for
 $\rb\in R^\max$ (i.e., with the open strata in the irreducible components of $C_X(V^\bullet)$).
 Similarly, $ \ol{\PP C}_X(V^\bullet) - \del\ol{\PP  C}_X(V^\bullet)$ is identified with
 the union of the strata $\PP C^\circ_{X,\rb}(V^\bullet)$ for
 $\rb\in R^\max$. 
 
 \end{enumerate}
  \end{thm}

  \vskip .2cm
 
 In fact, our construction provides, along the way, a wonderful compactification of
 any stratum in any   variety of complexes.
 
 \begin{thm}\label{thm:wonder2}
\begin{enumerate}
\item[(a)] The projection $C^{(|\rb|)}_{X,\rb}(V^\bullet)\to C_{X,\rb}(V^\bullet)$ is birational
 and biregular
over the open stratum  $C^\circ_{X,\rb}(V^\bullet)$. The subvarieties 
$C^{(|\rb|)}_{X,\rb}(V^\bullet)\cap C^{(|\rb|)}_{X,\rb'}(V^\bullet), \rb' < \rb$, are smooth and form
a divisor with normal crossings in $C^{(|\rb|)}_{X,\rb}(V^\bullet)$. 

\item[(b)] Similarly, 
the projection $\PP C^{(|\rb|)}_{X,\rb}(V^\bullet)\to \PP C_{X,\rb}(V^\bullet)$ is birational
 and biregular
over the open stratum  $\PP C^\circ_{X,\rb}(V^\bullet)$. The subvarieties 
$\PP C^{(|\rb|)}_{X,\rb}(V^\bullet)\cap \PP C^{(|\rb|)}_{X,\rb'}(V^\bullet), 0< \rb' < \rb$, are smooth and form
a divisor with normal crossings in $\PP C^{(|\rb|)}_{X,\rb}(V^\bullet)$. 

\end{enumerate}
 \end{thm}
 
 Our strategy for proving Theorems \ref{thm:blow} - \ref{thm:wonder2} consists in reducing
  the analysis of each blowup,
 by means of \'etale local charts, to the simplest case: the blowup of the zero section in the relative affine variety
 of complexes. This strategy  agrees with the general approach to  conical stratifications
 of singular varieties, sketched 
 in \S 3.3 of \cite{macpherson-procesi}. Our analysis also provides additional information that will be used in \S \ref{sec:CCSS}
  to identify $\k$-points 
 of the blowup with spectral sequences. 
 
 \vskip .2cm

 We start by analyzing this simplest case. 
 
 \vskip .2cm
 
 \paragraph  {D. Inductive step: structure of the first blowup.} 
 Consider $C_X^{(1)}(V^\bullet) = \Bl_X C_X(V^\bullet)$, where $X=C_{X,\0}(V^\bullet)$ is the zero section. 
 Since $C_X(V^\bullet)\to X$ is conic over $X$,    we have the projection $q: C_X^{(1)}(V^\bullet)\to \PP C_X(V^\bullet)$
 realizing $C_X^{(1)}(V^\bullet)$ as the total space of the relative line bundle $\Oc(-1)$. Therefore the strict transforms
 of the varieties $C_{X,\rb}(V^\bullet)$ (i.e., of the closures of strata) in $C_X^{(1)}(V^\bullet)$ are:
 \[
 C_{X, \rb}^{(1)}(V^\bullet) \,\,=\,\, q^{-1}\bigl( \PP C_{X, \rb}(V^\bullet)\bigr), \quad |\rb|\geq 1. 
 \]
 The lowest closures of the strata (coinciding with the corresponding strata) in $\PP C_X(V^\bullet)$
 are the $\PP C_{X,\rb}(V^\bullet)$ for $|\rb|=1$, i.e., for $\rb=e_i$ for some $i$. Being the strata, they are smooth and
 disjoint and carry vector bundles $H^\bullet_{\PP C_{X, \rb}(V^\bullet)}$, see \S  \ref{subsec:relcom}. 
 Let
 \[
 H^\bullet_\rb \,\,=\,\, q^*  \bigl(H^\bullet_{\PP C_{X, \rb}(V^\bullet)}\bigr), \quad |\rb|=1.
 \]
 This is a vector bundle on $C_{X,\rb}^{(1)}(V^\bullet)$. Proposition \ref{prop:chart} implies, now, by pullback:
 
 \begin{prop}\label{prop:chart2} 
 \noindent Let $|\rb|=1$. Then:
 \vskip .2cm
  
 (a)  Let $c$ be any $\k$-point of $C_{X,\rb}^{(1)}(V^\bullet)$. There exists an isomorphism of an
 \'etale neighborhood of $c$ in $C_{X,\rb}^{(1)}(V^\bullet)$ with an \'etale neighborhood of $c$ in
 the relative variety of complexes $C_{C_{X, \rb}^{(1)}(V)}(H^\bullet_\rb)$. Here, in the second variety,
 $c$ is understood as lying in the zero section.
 
 \vskip .2cm
 
 (b) Further, the  isomorphism in (a) can be chosen so that it takes, for any 
  $\rb'\geq \rb$,   an \'etale neighborhood of $c$ in $C^{(1)}_{X, \rb'}(V^\bullet)$
 to an \'etale neighborhood of $c$ in 
 $C_{ C_{X, \rb}^{(1)}  (V^\bullet ), \rb'-\rb }(H^\bullet_\rb)$. \qed
 \end{prop}
 
 \vskip .2cm
 
 \paragraph { E. Proof of Theorem \ref{thm:blow} and a ``disjointness lemma''.}
  We prove, by induction in $l$, the compound statement consisting
 of Theorem \ref{thm:blow} and the following claim.
 
 \begin{prop}\label{prop:chart3}
 Let $l=0,\cdots, d$. Let $|\rb|=l$ and $c$ be any $\k$-point of $C_{X, \rb}^{(l)}(V^\bullet)$ (resp. of 
 $\PP C_{X, \rb}^{(l)}(V^\bullet)$).
 
 \vskip .2cm 
 
 (a)  There exist:
 \begin{itemize}
 \item A Zariski open neighborhood $U$ of $c$ in $C_{X, \rb}^{(l)}(V^\bullet)$ (resp. in
 $\PP C_{X, \rb}^{(l)}(V^\bullet)$). 
 
 \item A graded vector bundle $H^\bullet$ on $U$.
 
 \item An isomorphism $\Xi=\Xi_{\rb, l}$ of  an \'etale neighborhood of $c$ in $C_X^{(l)}(V^\bullet)$
 (resp. in $\PP C_X^{(l)}(V^\bullet)$) with an \'etale neighborhood of $c$ in $C_U(H^\bullet)$. 
 
 \end{itemize}
 
 (b)  Further,  for any $\rb'\geq \rb$,  $\Xi$ restricts to an isomorphism of  an \'etale neighborhood of $c$ in
 $C_{X, \rb'}^{(l)}(V^\bullet)$ with  an \'etale neighborhood of $c$ in $C_{U, \rb'-\rb}(H^\bullet)$. 
 \end{prop}
 
 We assume the statements proven for a given value of $l$ (as well as for all the previous values). In particular,
 we define $C_X^{(l+1)}(V^\bullet)$ and $\PP C_X^{(l+1)}(V^\bullet)$
  by the formulas in Theorem \ref{thm:blow}(b).  After this we prove Theorem \ref{thm:blow} for $l+1$.
  
  \vskip .2cm
  
  The statement of Theorem  \ref{thm:blow}(a) for $l+1$ reads as follows.
  {\em
  For different $\rb$ with $|\rb|=l+1$, the subvarieties $C_{X, \rb}^{(l+1)}(V^\bullet)$, resp. 
  $\PP C_{X, \rb}^{(l+1)}(V^\bullet)$, are smooth and disjoint.} 
  
  To prove this, we note the following. By 
  Proposition \ref{prop:chart3}(b) 
  for $l$, we know that for any $\rb'$ with $|\rb'|=l$, 
  the variety $C_X^{(l)}(V^\bullet)$
  (resp.  $\PP C_X^{(l)}(V^\bullet)$) behaves near
   $C^{(l)}_{X, \rb'}(V^\bullet)$  (resp. near $\PP C^{(l)}_{X, \rb'}(V^\bullet)$),
  like the variety $C_U(H^\bullet)$ behaves near its zero section $U$ (\'etale local identification).
   So   $\Bl_{C^{(l)}_{X, \rb'}(V^\bullet)}C_X^{(l)}(V^\bullet)$
    (resp. $\Bl_{\PP C^{(l)}_{X, \rb'}(V^\bullet)}\PP C_X^{(l)}(V^\bullet)$)
 together with its closures of the strata, will be \'etale locally identified with $\Bl_U(C_U(H^\bullet))$
  together with its
 closures of the strata. This latter blowup was studied in Proposition \ref{prop:chart2}. In particular, 
 as pointed out in the discussion just above that proposition, its lowest
 closures of the strata, $C^{(1)}_{U, \rb''}(H^\bullet)$, $|\rb ''|=1$, are smooth and disjoint. But  Proposition
 \ref{prop:chart3}(b)  (for $l$)  implies that these lowest closures of the strata are \'etale locally identified
 with the strict transforms of the $C_{X, \rb}^{(l)}(V^\bullet)$, $ |\rb|=l+1$ in $C_X^{(l+1)}(V^\bullet)$, 
 i.e., with $C_{X, \rb}^{(l+1)}(V^\bullet)$.
 This proves part (a) of Theorem \ref{thm:blow} for $l+1$. 
 
 Part (b) of Theorem  \ref{thm:blow} for $l+1$ now constitutes the definition of $C_X^{(l+1)}(V^\bullet)$
 and $\PP C_X^{(l+1)}(V^\bullet)$.

 \vskip .2cm
 
  We now prove Proposition
 \ref{prop:chart3} for $l+1$. For this we  just need to combine two charts:
 \begin{enumerate}
 \item[(1)] The \'etale chart given by Proposition \ref{prop:chart2} for 
 $\Bl_U(C_U(H^\bullet))$
 
 \item[(2)]  The (source- and target-wise)  blowup of the already constructed 
  \'etale chart $\Xi_{\rb',l}$,  $|\rb'|=l$,  from Proposition  \ref{prop:chart3} for $l$.
 
 \end{enumerate}
 This concludes the inductive proof of Theorem \ref{thm:blow} and Proposition \ref{prop:chart3}. 
 \qed
 
 \vskip .2cm
 
 Let us conclude this part with the following ``disjointness lemma'', to be used later. 
 
 \begin{lem}\label{lem:disjointness}
 Let $\rb, \sb\in R$. Suppose that $\rb\not < \sb$ and $\sb\not < \rb$. Let $\rb'=\min(\rb,\sb)$
 (see Proposition \ref{prop:sch-strata-abs})
 and $l'=|\rb'|$. Then the varieties
 \[
 C^{(l'+1)}_{X, \rb}(V^\bullet), \,\,  C^{(l'+1)}_{X, \sb}(V^\bullet)\,\,\subset \,\, C_X^{(l'+1)}(V^\bullet)
 \]
 are disjoint. 
 \end{lem}
 
 \noindent{\sl Proof:} 
  By Proposition \ref{prop:C_X-properties}(d), $C_{X,\rb'}(V^\bullet) = C_{X,\rb}(V^\bullet)
 \cap C_{X,\sb}(V^\bullet)$ (scheme-theoretic intersection). Therefore the strict transforms
 of $C_{X,\rb}(V^\bullet)$ and $C_{X,\sb}(V^\bullet)$ in the blowup of $C_X(V^\bullet)$
 along $C_{X,\rb'}(V^\bullet)$ are disjoint by Proposition \ref{prop:blow-disjoint}. 
 Now, over the open stratum $C^\circ_{X, \rb'}(V^\bullet)$,
 this blowup coincides with $C_X^{(l'+1)}(V^\bullet)$ so we conclude that the  image of the intersection
 $C^{(l'+1)}_{X,\rb}(V^\bullet)\cap C^{(l'+1)}_{X,\sb}(V^\bullet)$ in $C_X(V^\bullet)$ does not meet
 $C^\circ_{X,\rb'}(V^\bullet)$.  
 
 \vskip .2cm
 
 It remains to eliminate the possibility of  a point  $p\in C^{(l'+1)}_{X,\rb'}(V^\bullet)$ belonging to
  $C^{(l'+1)}_{X,\rb}(V^\bullet)\cap C^{(l'+1)}_{X,\sb}(V^\bullet)$ and 
 projecting to a point    in  some smaller stratum $S$  inside  $C_{X,\rb'}(V^\bullet)$. Such
 a stratum has the form
  $S=C^\circ_{X,\tb}(V^\bullet)$ with $\tb < \rb'$.  
  
  \vskip .2cm

 Let $q$ be the image of $p$ in  $C^{(|t|)}_X(V^\bullet)$. 
 By  Proposition \ref{prop:chart3}, near $q$,
 the variety $C^{(|t|)}_X(V^\bullet)$ is \'etale locally identified with some $C_S(H^\bullet)$  so that
 $C_{X,\rb}^{(|t|)}(V^\bullet)$, resp. $C_{X,\sb}^{(|t|)}(V^\bullet)$, resp.
 $C_{X,\rb'}^{(|t|)}(V^\bullet)$
 is identified with 
 $C_{S,\rb-\tb}(H^\bullet)$, resp. $C_{S,\sb-\tb}(H^\bullet)$ resp. $C_{S,\rb'-\tb}(H^\bullet)$.
 Under this identification, the relevant part of
 $C_{X,\rb}^{(|t|+1)}(V^\bullet)$ is just the blowup of (the relevant part of) $C_{S,\rb-\tb}(H^\bullet)$
 along $C_{S,\rb'-\tb}(H^\bullet)$, and similarly for  $C_{X,\sb}^{(|t|+1)}(V^\bullet)$.
  Now observe that
 $\min(\rb-\tb, \sb-\tb) = \min(\rb, \sb)-\tb$. So, as before, the center of the blowup is the
 scheme-theoretic intersection of two subvarieties and so their strict transforms in the
 blowup are disjoint. We have thus shown that $C_{X,\rb}^{(|t|+1)}(V^\bullet)$ and $C_{X,\sb}^{(|t|+1)}(V^\bullet)$ do not intersect over $S$. Since $|t|<l'$, their subsequent iterated dominant transforms $C_{X,\rb}^{(l'+1)}(V^\bullet)$ and $C_{X,\sb}^{(l'+1)}(V^\bullet)$, respectively, do not intersect over $S$ as well. 
 
 \qed

 \paragraph{ F. Proof of Theorems \ref{thm:wonder1} and \ref{thm:wonder2}.} We start with some reductions.
 
 \vskip .2cm
 
 First, we will treat only the blowups $C_X^{(l)}(V^\bullet)$  of the affine varieties of complexes.
 The treatment of the  $\PP C_X^{(l)}(V^\bullet)$ is completely parallel. 
 
  \vskip .2cm
 
 Second, note that Theorem \ref{thm:wonder1} is a particular case of Theorem \ref{thm:wonder2}.
 Indeed, each $C^{(d)}_{X,\rb}(V^\bullet)$, $\rb\in R^\max$, will first appear as $C_{X, \rb}^{(|\rb|)}(V^\bullet)$
 and will not change in the subsequent blowups. So we concentrate on the proof of 
  Theorem \ref{thm:wonder2}.
  
   \vskip .2cm
  
  Let us write $C^{(l)}_\rb = C^{(l)}_{X, \rb}(V^\bullet)$. For any $\rb\in R$ and any $l\leq |\rb|$ put
  \be\label{eq:def-W}
  W_\rb^{(l)} \,\,=\,\, C_\rb^{(l)} \,\,-\,\,\bigcup_{\sb <\rb, \,\,\, l\leq |\sb|} \bigl(C_\rb^{(l)}\cap C_\sb^{(l)}\bigr). 
  \ee
 This is an open subvariety in $C_\rb^{(l)}$. Put also
 \[
 D_\rb^{(l)} \,\,=\,\, W_\rb^{(l)} \,\,\cap \bigcup_{\sb <\rb, \,\,\, l> |s|} C_\sb^{(l)}.
 \]
 When $l=0$, we have  that $W_\rb^{(0)}=C_\rb^\circ = C^\circ_{X,\rb}(V^\bullet)$ is the open stratum
 corresponding to $\rb$
 in the variety of complexes, while $D_\rb^{(0)}=\emptyset$. 
 
 \vskip .2cm
 
 When $l=|\rb|$, we have that $W_\rb^{(|\rb|)} = C_\rb^{(|\rb|)}$, and 
 \[
 D_\rb^{(|\rb|)} \,\,=\,\, C_\rb^{(|\rb|)} \,\,\cap \,\,\bigcup_{\sb< \rb} C_\sb^{(|\rb|)}
 \]
 is the divisor which is claimed in the theorem to be a divisor with normal crossings.
 So it suffices to prove the following more general statement. In this statement and its
 proof we will use the following terminology. A pair $(Z,D)$ will be called {\em wonderful}, if $Z$
 is a smooth variety and $D$ is a divisor with normal crossings in $Z$. 
 
 \begin{prop}\label{prop:DNC}
 For any $\rb\in R$ and any $l\leq |\rb|$, the pair $(W^{(l)}_\rb, D^{(l)}_\rb)$ is wonderful. 
 \end{prop}
 
 \noindent {\sl Proof:} We proceed by induction in $l$. The case $l=0$ is clear from the above.
 Suppose the statement is proved for a given value of $l$, and suppose that $(l+1)\leq |\rb|$,
 so that the next statement is a part of the proposition. 
 Look at the blowup map $p: C^{(l+1)}\to C^{(l)}$ with the smooth center
 $\coprod_{|\sb|=l} C_\sb^{(l)}$, as in Theorem \ref{thm:blow}(b). 
 
 \begin{lem}
 $p$ is biregular over $W^{(l)}_\rb\subset C^{(l)}$, i.e., $W^{(l)}_\rb$ does not meet the center of
 the blowup.
 \end{lem}
 
 \noindent {\sl Proof of the lemma:} Indeed, each $C_\sb^{(l)}$, $|\sb|=l$, will either not
 meet $C^{(l)}_\rb$ and hence $W^{(l)}_\rb$ (this will happen if $\sb\not < \rb$ by Lemma
 \ref{lem:disjointness}), or will
 meet $C_\rb^{(l)}$ but will be removed in forming $W_\rb^{(l)}$ (this will happen if $\sb < \rb$). 
 \qed
 
 \vskip .2cm
 
 Denote by $E$ the exceptional divisor of $p$ (the preimage of the center of the blowup).
 The lemma means that any  ``new" point   $w\in W^{(l+1)}_\rb$ (i.e., a point not lifted by a
 local biregular map from a point in $W_\rb^{(l)}$), lies in $E$. So it is enough to prove that
  $(W^{(l+1)}_\rb, D^{(l+1)}_\rb)$ is wonderful only near such new points $w$, belonging to $E$.

  \vskip .2cm
  
  So we choose such $w$ and denote $c=p(w)$. Then $c\in C_\sb^{(l)}$ for some $\sb$
  with $|\sb|=l$. Since $w\in W_\rb^{(l+1)}$, we have $\sb < \rb$. We now apply Proposition
  \ref{prop:chart3} to get an open neighborhood $U$ of $c$ in $C_\sb^{(l)}$, a graded vector
  bundle $H^\bullet$ on $U$ and 
  an identification $\Xi$  of an \'etale neighborhood of $c$ in $C_\rb^{(l)}$
  with an \'etale neighborhood of $c$ in $C_{U, \rb-\sb}(H^\bullet)$. Applying the blowup
  along the intersection of $U$ with the \'etale neighborhoods in the source and target of $\Xi$,
  we identify an \'etale neighborhood of $w$ in $C_\rb^{(l+1)}$ with the \'etale neighborhood 
  of a point $w'$ in $\Bl_U(C_{U, \rb-\sb}(H^\bullet))$, which is the total space of the line
  bundle $\Oc(-1)\to \PP C_{U, \rb-\sb}(H^\bullet)$.
  
  \begin{lem}
  (a) The point $w'$ lies on the zero section of the bundle $\Oc(-1)$. 
  
  \vskip .2cm
  
  (b) Identifying this zero section with $\PP C_{U, \rb-\sb}(H^\bullet)$, we have that $w'$
  lies in the open stratum $\PP C^\circ_{U, \rb-\sb}(H^\bullet)$. 
  \end{lem}
  
  \noindent {\sl Proof of the lemma:}  (a) follows because $w$ lies in the exceptional divisor of $p$.
  
  \vskip .2cm
  
  (b) Applying  \eqref{eq:def-W} in our case, we can write that
  \[
  w\in W_\rb^{(l+1)} \,\,=\,\, C_\rb^{(l+1)} - \bigcup_{\rb'< \rb \atop
  l+1\leq |\rb'| } (C_\rb^{(l+1)} \cap C_{\rb'}^{(l+1)})\,\,\subset \,\, C_\rb^{(l+1)}. 
  \]
  Since $w\in E$, it represents a point in the projectivization of the normal cone
  \[
  \NC_{C_\sb^{(l)}} W_\rb^{(l)} \,\,\subset\,\, 
  \NC_{C_\sb^{(l)}} C_\rb^{(l)}.
  \]
  The variety $\NC_{C_\sb^{(l)}} C_\rb^{(l)}$ is identified under $\Xi$ (\'etale locally around $c\in C_\rb^{(l)}$) with $C_{\rb-\sb}(H^\bullet)$. 
  Under this identification, the parts removed in forming $W_\rb^{(l+1)}$, namely $C_\rb^{(l+1)} \cap C_{\rb'}^{(l+1)}$
  match the subvarieties $\PP C_{\rb'-\sb}(H^\bullet)$. More precisely, the intersection 
  $C_\rb^{(l+1)} \cap C_{\rb'}^{(l+1)}\cap E$, is identified with $\PP C_{\rb'-\sb}(H^\bullet)$.
  So if the statement of part (b) is not true, then $w$ would lie in one of the removed parts. \qed

  \vskip .2cm

 \noindent Now Proposition \ref{prop:DNC} follows from the next obvious statement.
 
 \begin{lem}
 Let $(Z,D)$ be a wonderful pair, and $q: L\to Z$ be a line bundle.
 Then $(L, q^{-1}(D)\cup Z)$ (with $Z\subset L$ being the zero section), is a wonderful pair. \qed
 \end{lem}
 
\noindent Therefore, Theorems \ref{thm:wonder1} and \ref{thm:wonder2} are proved.


 \section{Complete complexes and spectral sequences}\label{sec:CCSS}
 
 \paragraph{A. Single-graded spectral sequences.}
 By a {\em spectral sequence} of $\k$-vector spaces we mean a sequence of complexes
 $(E^\bullet_\nu, D^\nu)\in\Com_\k$, $\nu=0,\cdots, k+1$ such that
 $E_{\nu+1}^{\bullet} = H^\bullet_{D^\nu}(E_\nu^\bullet)$ for each $\nu <k+1$. Here $k+1$ can be either a finite
 number of $\infty$.  If $k+1$ is finite, then the differential $D^{k+1}$ is considered to be zero (so that there is no
 additional $E_{k+2}^\bullet$ to speak of).

\begin{Defi}\label{def:SS-red}
 A spectral sequence $(E^\bullet_\nu, D^\nu)$, $\nu=0,\cdots, k+1$,  will be called {reduced}, if:
 \begin{enumerate}
 \item[(1)]  $E_0^\bullet$ (and therefore each $E_\nu^\bullet$)
 is  finite-dimensional, i.e., the total dimension $\dim(E_0^\bullet)<\infty$.
 
 \item[(2)] Each $D^\nu$, $\nu=1, \cdots, k$,  is not entirely $0$,
 i.e., at least one component $D^\nu_i: E_\nu^i\to E_\nu^{i+1}$ is nonzero. 
 
 \item[(3)] The graded vector space $E_{k+1}^\bullet$ is sparse, see \S \ref{sec:blowups}B. 
 \end{enumerate}
 We say that $(E^\bullet_\nu, D^\nu)$ is  strongly reduced if, in addition,  $D^0\neq 0$. 
 
 \end{Defi}
 For a reduced spectral sequence we have  $\dim (E_{\nu+1}^\bullet) < \dim(E_\nu^\bullet)$, 
 so the length of such a sequence is bounded by $\dim(E_0^\bullet)$.

 \begin{Defi}
  Let $V^\bullet$ be a finite-dimensional graded $\k$-vector space. A {complete complex} (affine version) 
  on $V^\bullet$ is an equivalence class of   reduced spectral sequences $(E_\nu^\bullet, D^\nu)$ with $E_0^\bullet=V^\bullet$,
 where  each differential $D^\nu, \nu\geq 1$ is considered modulo  scaling (the same scalar for all components
 $D^\nu_i$). 
 
  A {complete complex} (projective version) 
  on $V^\bullet$ is an equivalence class of   strongly reduced spectral sequences 
  $(E_\nu^\bullet, D^\nu)$ with $E_0^\bullet=V^\bullet$,
 where  each differential $D^\nu, \nu\geq 0$ is considered modulo  scaling (the same scalar for all components
 $D^\nu_i$).

 \end{Defi}

 We denote by  $\SS(V^\bullet)$ and $\PP \SS(V^\bullet)$ the sets of complete complexes on $V^\bullet$ 
 in the affine and projective version respectively.

 \paragraph{B. $\k$-points of $\ol C(V^\bullet)$ and $\ol{\PP C}(V^\bullet)$ as spectral sequences.}
 Let $C^\circ(V^\bullet)$ be the generic part of $C(V^\bullet)$, i.e., the union of the maximal strata,
 see \eqref{eq:generic-str}. 
 We have an embeding  $C^\circ(V^\bullet)(\k)\subset \SS(V^\bullet)$:
 a differential  $D$ making $V^\bullet$ into a complex, is identified
 with a spectral sequence of length 1, that is, consisting only of $E_0^\bullet = V^\bullet$
 and $E_1^\bullet = H^\bullet_{D}(E_0^\bullet)$. The fact that $D$ lies in a maximal stratum
 means,  by Proposition \ref{prop:max-sparse}, means that $E_1^\bullet$ is sparse, so
 condition (3) of Definition \ref{def:SS-red} is satisfied. 
  We have a similar embedding
 $\PP C(V^\bullet)(\k)\subset \PSS(V^\bullet)$.
 
 \begin{thm}\label{thm:CCSS}
 We have  identifications 
 \[
 \ol C(V^\bullet)(\k) \simeq \SS(V^\bullet), \quad \ol{\PP C}(V^\bullet)(\k)
 \simeq \PSS(V^\bullet), 
 \]
  extending the above embeddings. 
 \end{thm}
 
 \paragraph{C. Stratification of complete complexes.} 
 Before proving the theorem, we study the natural stratifications of
 $\ol C(V^\bullet)$ and $\ol {\PP C}(V^\bullet)$ given by the generic parts of all possible intersections of the boundary divisors in each of these wonderful compactifications. It is convenient to work in the relative situation
 of the relative varieties of complete complexes $\ol C_X(V^\bullet)$ and $\ol {\PP C}_X(V^\bullet)$
 corresponding to a graded vector bundle $V^\bullet$ on a smooth
 variety $X$. We recall the divisors $\Delta_\rb\subset \ol C_X(V^\bullet)$, $\rb\in R^\flat$ and $\PP\Delta_\rb\subset \ol {\PP C}_X(V^\bullet)$, $\rb\in R^\flat_{>0}$ from 
 Theorem \ref{thm:wonder1}. To emphasize their dependence on $X, V^\bullet$
 we will write $\Delta_\rb^{\ol C_X(V^\bullet)}$ and $\PP\Delta_\rb^{\ol {\PP C}_X(V^\bullet)}$ respectively, if needed.

 \begin{prop}
 Let $\rb^{(1)}, \cdots, \rb^{(k)}\in R^\flat$ (resp. $\rb^{(1)}, \cdots, \rb^{(k)}\in R^\flat_{>0}$) be distinct. 
 The intersection   $\Delta_{\rb^{(1)}}\cap \Delta_{\rb^{(2)}}\cap\cdots \cap \Delta_{\rb^{(k)}}\subset \ol C_X(V^\bullet)$ (resp. $\PP\Delta_{\rb^{(1)}}\cap \PP\Delta_{\rb^{(2)}}\cap\cdots \cap \PP\Delta_{\rb^{(k)}}\subset \ol {\PP C}_X(V^\bullet)$)
 is nonempty if and only if, after a permutation of the $\rb^{(i)}$, we have
 $\rb^{(1)} < \cdots < \rb^{(k)}$. 
 \end{prop}
 
 \noindent{\sl Proof:} We only prove the statement about the  variety of complete complexes; the projective case follows by identical arguments. \\
 \indent``If":  We proceed by induction on $k$, the case $k=1$ being trivial.
 So we assume the statement proved for all $X, V^\bullet$ and $\rb^{(1)} <\cdots
 < \rb^{(k-1)}$. 
 
 Suppose  now some $X, V^\bullet$ and 
  $\rb^{(1)} < \cdots < \rb^{(k)}$ are given. We consider the stratum $S= C^\circ_{X,\rb^{(1)}}(V)$.
  By Propositions \ref{prop:rel-ncone} and \ref{prop:chart},  each point of $S$ has an \'etale neighborhood $U \to C_X(V^\bullet)$ 
  identified with a part (\'etale) of $C_S(H^\bullet)$ where $H^\bullet$ is the 
  vector bundle of the cohomology on $S$. 
  Under this identification, each subvariety $C_{S, \sb}(H^\bullet)$ corresponds to
  $C_{X, {\rb^{(1)}+\sb}}(V^\bullet)$. 
  
  \vskip .2cm

  Accordingly, the preimage of $U$ in $\ol C_X(V^\bullet)$ is identified with a part of
 $\ol C_S(H^\bullet)$ in such a way that the divisors $\Delta_\sb^{\ol C_S(H^\bullet)}$ in
  $\ol C_S(H^\bullet)$ correspond to
  the divisors $\Delta^{\ol C_X(V^\bullet)}_{\rb^{(1)}+\sb}$ in $\ol C_X(V^\bullet)$. 
  In particular, $\Delta_{\rb^{(1)}}^{\ol C_X(V^\bullet)}$ corresponds to the dominant transform
  of the zero section of $C_S(H^\bullet)$ which is nothing but $\ol{\PP C}_S(H^\bullet)$,
  the projective variety of complete complexes.   
  
  \vskip .2cm
  
  Since, by the inductive assumption,
  the intersection $\Delta_{\rb^{(2)}-\rb^{(1)}} \cap \cdots \cap \Delta_{\rb^{(k)}-\rb^{(1)}}$
  in $\ol C_S(H^\bullet)$ is nonempty, the intersection of their images in
  $\ol{\PP C}_S(H^\bullet)$ is also non-empty. But by the above argument, this
  intersection in  $\ol{\PP C}_S(H^\bullet)$ is \'etale locally identified with a part of 
    the intersection
 $\Delta_{\rb^{(1)}}\cap \Delta_{\rb^{(2)}}\cap\cdots \cap\Delta_{\rb^{(k)}}$ in $\ol C_X(V^\bullet)$,
  which is therefore nonempty too. 
  
  \vskip .3cm
  
  ``Only if": 
 The statement reduces to the following:  if $\Delta_\rb\cap \Delta_\sb\neq\emptyset$, then
  $\rb < \sb$ or $\sb < \rb$.  To prove this, suppose that $\rb\not < \sb$ and $\sb\not< \rb$. 
Let $\rb' =  \min(\rb, \sb)$, see Proposition \ref{prop:sch-strata-abs}. 
 and   $l= |\rb'|$. 
 Then $\rb' < \rb, \sb$. Our statement now follows from Lemma \ref{lem:disjointness}. 
  \qed
 
 \vskip .2cm
 We now make precise the natural stratification of the varieties of complete complexes associated to their boundary. To this end, we give the following Definition:
 \begin{Defi}\label{defi:complete-strata}
Let $T\subset R$ be any subset of the form $T=\{\rb^{(1)} < \cdots < \rb^{(k)}\}$, where $\rb^{(1)}, \cdots, \rb^{(k)}\in
 R^\flat$. 
We allow the case $k=0$, i.e., $T=\emptyset$. 
\begin{enumerate}\item The stratum of $\ol C_X(V^\bullet)$ associated to $T$ is defined to be the locally closed subvariety 
$$
\Delta^\circ_T (X, V^\bullet) \,\, = \,\, \Delta_{\rb^{(1)}, \cdots, \rb^{(k)}} (X, V^\bullet) \,\,:=\,\,
\bigcap_{\rb\in T} \Delta_\rb \,\,- \,\,\bigcup_{\sb\notin T} \Delta_\sb \,\,\subset \,\, \ol C_X(V^\bullet).
$$
 \item Let $\rb^{(1)}\neq {\bf 0}$.  The stratum of $\ol {\PP C}_X(V^\bullet)$ associated to $T$ is 
 defined to be the locally closed subvariety 
$$
\PP \Delta^\circ_T (X, V^\bullet) \,\, = \,\, \PP\Delta_{\rb^{(1)}, \cdots, \rb^{(k)}} (X, V^\bullet) \,\,:=\,\,
\bigcap_{\rb\in T}\PP \Delta_\rb \,\,- \,\,\bigcup_{\sb\notin T} \PP\Delta_\sb \,\,\subset \,\, \ol {\PP C}_X(V^\bullet).
$$
\end{enumerate}
  \end{Defi}
 
 \begin{rem} If $T=\emptyset$, the stratum of $\ol C_X(V^\bullet)$ (resp. $\ol {\PP C}_X(V^\bullet)$) associated to $T$ is the generic part $C_X^\circ(V^\bullet)$ (resp. $\PP C_X^\circ(V^\bullet)$) of $\ol C_X(V^\bullet)$ (resp. $\ol {\PP C}_X(V^\bullet)$) that is, the union of $C^\circ_{X,\rb}(V^\bullet)$ (resp. $\PP C^\circ_{X,\rb}(V^\bullet)$) for all $\rb\in R^\max$.
\end{rem}

 If $X=\Spec(\k)$, i.e., $V^\bullet$ is just a graded $\k$-vector space, we  abbreviate the notation for
 the above varieties to  $\Delta^\circ_T(V^\bullet)$,
 resp. $\PP\Delta^\circ_T(V^\bullet)$.
 
 \paragraph{D. Proof of Theorem \ref{thm:CCSS}.}
 Let $V^\bullet$ be a graded $\k$-vector space, as in the theorem. 
 We will identify each stratum in  $\ol C (V^\bullet)$, resp. $\ol{\PP C} (V^\bullet)$, with the set of
 spectral sequences with fixed numerical invariants. Let $T=\{\rb^{(1)} < \cdots < \rb^{(k)}\}$ be as above.
 Define the set $SS^\circ_T(V^\bullet)\subset \SS(V^\bullet)$ to consist of  equivalence classes of 
 spectral sequences $(E_\nu^\bullet, D^\nu)$ (see discussion after Definition \ref {def:SS-red}) such that:
 \begin{enumerate}
 \item[(0)] $E_0^\bullet = V^\bullet$ and $D^0\in C^\circ_{\rb^{(1)}}(E_0^\bullet)$;
 
 \item[(1)] $D^1\in C^\circ_{\rb^{(2)}-\rb^{(1)}}(E_1)$, where $E^\bullet_1 :=H^\bullet_{D^0}(E_0^\bullet)$, 
 
 \item[(2)] $D^2\in C^\circ_{\rb^{(3)}-\rb^{(2)}}(E_2)$, and so on. 
 \end{enumerate}
 If $\rb^{(1)}\neq {\bf 0}$, we  denote by $\PP\SS^\circ_T(V^\bullet)$ the subset of $\PP\SS(V^\bullet)$
 corresponding to $SS^\circ_T(V^\bullet)$.  It is clear that we have disjoint decompositions
 \[
 \SS(V^\bullet) \,\,=\,\,\bigsqcup_{T=\{\rb^{(1)} < \cdots < \rb^{(k)}\}} \SS^\circ_T(V^\bullet), 
 \quad \PP\SS(V^\bullet)   \,\,=\,\,\bigsqcup_{T=\{0< \rb^{(1)} < \cdots < \rb^{(k)}\}} \PP\SS^\circ_T(V^\bullet). 
 \]
 Theorem \ref{thm:CCSS} is a consequence of the following refined statement.
 
 \begin{prop}\label{prop:strata=SS}
 For any $T=\{\rb^{(1)} < \cdots < \rb^{(k)}\}$ as above we have  identifications
 \[
 \SS^\circ_T(V^\bullet) \,\,\simeq \,\, \Delta^\circ_T(V^\bullet) (\k) , \quad \PP \SS^\circ_T(V^\bullet) \,\,\simeq \,\, \PP\Delta^\circ_T(V^\bullet) (\k).
 \]
  \end{prop}
 
 \noindent{\sl Proof of the proposition:} For the proof, we work with relative  complete varieties of complexes
 and analyze their strata, introduced in part C, in an inductive fashion. 
 
 \begin{lem}\label{lem:SS-strata}
 Let $X$ be a smooth  variety over $\k$ and $V^\bullet$ a graded vector bundle on $X$. Then we have 
 an isomorphism
 \[
 \Delta^\circ_{ \rb^{(1)} , \cdots , \rb^{(k)}}(X, V^\bullet) \,\,\simeq \,\,
 \PP\Delta^\circ_{ \rb^{(2)}-\rb^{(1)}, \rb^{(3)}-\rb^{(1)}, \cdots ,  \rb^{(k)}-\rb^{(1)} }  (C^\circ_{X,\rb^{(1)}} (V^\bullet),
 H^\bullet),
 \]
 where $H^\bullet$ is the vector bundle of cohomology on the stratum $C^\circ_{X,\rb^{(1)}} (V^\bullet)$.
 We further have an isomorphism
 \[
 \PP\Delta^\circ_{ \rb^{(1)} ,\cdots , \rb^{(k)}}(X, V^\bullet) \,\,\simeq \,\,
 \PP\Delta^\circ_{ \rb^{(2)}-\rb^{(1)}, \rb^{(3)}-\rb^{(1)}, \cdots ,  \rb^{(k)}-\rb^{(1)} }  (\PP C^\circ_{X,\rb^{(1)}}  (V^\bullet),
 H^\bullet).
 \]
  \end{lem}

  Knowing the lemma, the proof of Proposition \ref{prop:strata=SS} (and thus of Theorem  \ref{thm:CCSS})
  for  strata in $\ol C_X(V^\bullet)$ 
  is finished as follows. 
  We construct inductively the following varieties $X_\nu$ together with graded vector bundles $E_\nu^\bullet$
  on them:
  \begin{enumerate}
  \item[(0)] $X_0=\Spec(\k)$,  and $E_0^\bullet = V^\bullet$.
  
  \item[(1)] $X_1$ is the stratum $\PP C^\circ_{\rb^{(1)}, X_0}(E_0^\bullet)$, and $E_1^\bullet$ is the bundle of 
  cohomology on this stratum.
  
  \item[(2)] $X_2$  is the stratum $\PP C^\circ_{\rb^{(2)}-\rb^{(1)}, X_1}(E_1^\bullet)$, and $E_2^\bullet$ is the bundle of 
  cohomology on this stratum,  
  
  \item[] .......
    
  \item[($k$)] $X_k$  is the stratum $\PP C^\circ_{\rb^{(k)}-\rb^{(k-1)}, X_{k-1}}(E_{k-1}^\bullet)$, and $E_k^\bullet$ is the bundle of  cohomology on this stratum. 
    \end{enumerate}
 
 Lemma \ref{lem:SS-strata} implies, by induction, the following:
 
 \begin{cor}
 The variety $X_k$ is identified with $\Delta^\circ_T(V^\bullet)$. \qed
 \end{cor}
 
 Proposition \ref{prop:strata=SS}  for  strata in $\ol C_X(V^\bullet)$ 
  now follows because points of $X_k$ are manifestly identified with
 equivalence classes of spectral sequences, as in Definition \ref{def:SS-red}. 
 The case of strata in $\ol{\PP C}_X(V^\bullet)$ is treated similarly.

\paragraph {E.  Proof of Lemma   \ref{lem:SS-strata}.}  By definition,
 $\Delta_{\rb^{(1)}} = \Delta_{\rb^{(1)}}(X, V^\bullet)$ is the iterated dominant transform of the closed subvariety
 $C_{X,\rb^{(1)}}(V^\bullet)$ in the first tower of blowups in Theorem \ref{thm:blow}. 
 It follows that
 \[
 \Delta_{\rb^{(1)}} - \bigcup_{\sb < \rb^{(1)}} \Delta_{\sb} \,\,=\,\, \wt{C^\circ}_{X,\rb^{(1)}}(V^\bullet)
 \]
is the iterared dominant transform of the open part (stratum) $C^\circ_{X,\rb^{(1)}}(V^\bullet)$. Let us analyze
this iterated transform and the tower of blowups in more detail. 

The  first  blowup   nontrivial over  $C^\circ_{X,\rb^{(1)}}(V^\bullet)$,  will appear at the stage $l=|\rb^{(1)}|$.
 It will
be the blowup
along the dominant transform of  $C_{X,\rb^{(1)}}(V^\bullet)$ which, on our part, reduces to
 $C^\circ_{X,\rb^{(1)}}(V^\bullet)$ itself. 
The corresponding dominant transform is, therefore,  the total inverse image, i.e., the projectivization of the
normal cone to  $C^\circ_{X,\rb^{(1)}}(V^\bullet)$ in $C_X(V^\bullet)$. This projectivization is the 
projective variety of complexes
$\PP C_{C^\circ_{X,\rb^{(1)}}(V^\bullet)} (H^\bullet)$. 

If we continue the construction of $\ol C_X(V^\bullet)$
in the tower of blowups of Theorem \ref{thm:blow} ,
then subsequent  blowups along dominant transforms
of the $C_{X,\tb}(V^\bullet)$, $\tb > \rb^{(1)}$ will induce blowups of 
$\PP C_{C^\circ_{X,\rb^{(1)}}(V^\bullet)} (H^\bullet)$ along (dominant transforms of) the 
$\PP C_{C^\circ_{X,\rb^{(1)}}(V^\bullet),\tb-\rb^{(1)}} (H^\bullet)$. This will produce 
$\ol{\PP C}_{C^\circ_{X,\rb^{(1)}}(V^\bullet)}(H^\bullet)$, 
the relative projective variety of  complete complexes. In other words, we have established an identification
\[
 \Delta_{\rb^{(1)}} - \bigcup_{\sb < \rb^{(1)}} \Delta_{\sb} \,\,\simeq 
 \,\, \ol{\PP C}_{C^\circ_{\rb^{(1)}, X}(V^\bullet)}(H^\bullet). 
\]
Under this identification the intersection of the LHS  with each  divisor $\Delta_\tb^{{\ol C_X(V^\bullet)}}$, $\tb >\rb^{(1)}$,  
corresponds to  the divisor $\PP\Delta_{\tb-\rb^{(1)}}$ in  $\ol{\PP C}_{C^\circ_{\rb^{(1)}, X}(V^\bullet)}(H^\bullet)$. 
The lemma is immediate from this. \qed


\vskip 1cm

\noindent{Kavli IPMU (WPI), UTIAS, The University of Tokyo, Kashiwa, Chiba 277-8583, Japan.} \\
Email: {\tt mikhail.kapranov@ipmu.jp}  and {\tt evangelos.routis@ipmu.jp}

\ed

\end{document}